\documentclass[opre,nonblindrev]{oo}

\DoubleSpacedXI

\usepackage{endnotes}
\let\footnote=\endnote

\usepackage{natbib}
\bibpunct[, ]{(}{)}{,}{a}{}{,}%
\def\bibfont{\small}%

\usepackage[toc,page]{appendix}

\TheoremsNumberedThrough     
\ECRepeatTheorems
\EquationsNumberedThrough    

\usepackage{amsmath,amssymb,amsfonts}  
\usepackage[dvipsnames]{xcolor}
\usepackage{tikz}
\usetikzlibrary{arrows,shapes,fit}
\usetikzlibrary{fadings}
\usetikzlibrary{patterns}
\usetikzlibrary{shadows.blur}
\usepackage{extarrows}
\usepackage{physics}
\usepackage{mathdots}
\usepackage{yhmath}
\usepackage{cancel}
\usepackage{color}
\usepackage{siunitx}
\usepackage{array}
\usepackage{gensymb}
\usepackage{tabularx}

\usepackage{booktabs} 
\usepackage{comment}
\usepackage{url}
\usepackage{eurosym}
\usepackage{graphicx}
\usepackage{multicol}
\usepackage{multirow}
\usepackage{floatrow}
\usepackage{float}
\usepackage{caption}
\usepackage{subcaption}
\usepackage{rotating}
\usepackage{algorithm}
\usepackage[compatible]{algpseudocode}
\usepackage{natbib}
\usepackage[colorlinks=true,linkcolor=blue,citecolor=cyan, hypertexnames=false]{hyperref}
\usepackage{lineno}
\modulolinenumbers[5]

\allowdisplaybreaks

\newcommand{\convx}{\textup{conv}}
\newcommand{\conv}[1]{\convx (#1)}
\newcommand{\dia}[1]{\textup{dia} (#1)}

\newcommand{\Exp}{\mathbb{E}}
\newcommand{\R}{\mathbb{R}}
\newcommand{\Z}{\mathbb{Z}}

\newcommand{\dist}{\mathbb{P}}
\newcommand{\Kpar}{K}
\newcommand{\Kt}{\Kpar_t}
\newcommand{\Ksetprime}{K'}
\newcommand{\Kprimet}{{\Ksetprime_t}}
\newcommand{\kprimeindex}{{k'}}
\newcommand{\ntminusone}{\npar_{t-1}}
\newcommand{\mpar}{m}
\newcommand{\mstatet}{\mpar^\textsc{s}_t}
\newcommand{\mrecourset}{\mpar^\textsc{r}_t}
\newcommand{\npar}{n}
\newcommand{\nt}{\npar_t}

\newcommand{\xipar}{\xi}
\newcommand{\xiparrand}{\boldsymbol{\xipar}}
\newcommand{\xiparval}{\xipar}
\newcommand{\xisubtrand}{\xiparrand_t}
\newcommand{\xisubtval}{\xiparval_t}

\newcommand{\xicolor}{black} 
\newcommand{\xiprimecolor}{black} 
\newcommand{\xiaveragecolor}{\color{black}} 
\newcommand{\xitwostagecolor}{\color{black}} 

\newcommand{\xiT}{{\color{\xicolor}{\xiparrand^T}}}

\newcommand{\xiprimeT}{{\color{\xiprimecolor}{{\xiparrand'}^T}}}
\newcommand{\xisupt}{{\color{\xicolor}{\xiparrand^t}}}
\newcommand{\xisups}{{\color{\xicolor}{\xiparrand^s}}}
\newcommand{\xisupsminusone}{{\color{\xicolor}{\xiparrand^{s-1}}}}
\newcommand{\xiprimesupt}{{\color{\xiprimecolor}{{\xiparrand'}^t}}}
\newcommand{\xisuptminusone}{{\color{\xicolor}{\xiparrand^{t-1}}}}
\newcommand{\xisuptplusone}{{\color{\xicolor}{\xiparrand^{t+1}}}}
\newcommand{\xiprimesuptminusone}{{\color{\xiprimecolor}{{\xiparrand'}^{t-1}}}}
\newcommand{\xiprimesuptplusone}{{\color{\xiprimecolor}{{\xiparrand'}^{t+1}}}}

\newcommand{\xiTval}{{\color{\xicolor}{\xiparval^T}}}

\newcommand{\xisuptval}{{\color{\xicolor}{\xiparval^t}}}
\newcommand{\xisupomegaval}{{\color{\xicolor}{\xiparval^T_\omega}}}
\newcommand{\xisuptomegaval}{{\color{\xicolor}{\xiparval^t_\omega}}}
\newcommand{\xisuptmoneomegaval}{{\color{\xicolor}{\xiparval^{t-1}_\omega}}}

\newcommand{\xisuptminusoneval}{{\color{\xicolor}{\xiparval^{t-1}}}}

\newcommand{\xiaverageval}{{\xiaveragecolor{\bar{\xiparval}^s_{|\xisupt}}}}
\newcommand{\xitwostagerand}{{\xitwostagecolor{{\xiparrand}_{|t}}}}

\newcommand{\xitwostagevalomega}{{\xitwostagecolor{{\xiparval}_{|t,\omega}}}}

\newcommand{\allxi}{\mathbb{P}\text{-a.e.} \ \xiT\in \Xi}
\newcommand{\allxiscen}{\mathbb{P}\text{-a.e.} \ \xiTval\in \Xi}

\newcommand{\MSLP}{MSLP}
\newcommand{\MSMIP}{MSMIP}

\newcommand{\ldr}{\beta}
\newcommand{\ldrvec}{\ldr}
\newcommand{\ldrval}{\ldr}
\newcommand{\naldr}{\alpha}
\newcommand{\naldrvec}{\naldr}
\newcommand{\naldrval}{\naldr}

\newcommand{\set}[1] {[#1]}
\newcommand{\bsub}{\begin{subequations}}
\newcommand{\esub}{\end{subequations}}

\newcommand{\spacing}{\DoubleSpacedXI}

\newcommand{\rev}[1]{\noindent{#1}}
\newcommand{\revtwo}[1]{\noindent{#1}}

\newcommand{\bsubeq}{\begin{subequations}}
\newcommand{\esubeq}{\end{subequations}}
\newcommand{\BI}{\begin{itemize}}
\newcommand{\EI}{\end{itemize}}
\newcommand{\BE}{\begin{enumerate}}
\newcommand{\EE}{\end{enumerate}}
\newcommand{\I}{\item}
\newcommand{\BSE}{\begin{subequations}}
\newcommand{\ESE}{\end{subequations}}

\newcommand{\Graph}{G}
\newcommand{\nodeSet}{\mathcal{V}}
\newcommand{\node}{v}
\newcommand{\edgeSet}{\mathcal{L}}
\newcommand{\edge}{\ell}
\newcommand{\waveSet}{\Omega}
\newcommand{\wavelength}{\omega}

\def \SD {\mathcal{SD}}
\def \SDcurr {\SD^1} 

\def \abilene {ABILENE}

\newcommand{\edgeSetw}{\mathcal{L}(\wavelength)}
\def \DnewMSIP {r^{t,(s,d)}(\xisupt)}
\def \DnewMSIPBF {r^{t,(s,d)}}

\def \SDsceMSIP {\mathcal{SD}^t(\xisupt)}
\def \SDsceMSIPBF {\SD^t}

\def \SDlsceMSIPBF {\SD^t_{\hat{\edge}}}

\def \linkMSIP {x_{\edge,\wavelength}^{1,(s,d)}(\xiparrand^1)}
\def \wavelinkMSIP {\mu_{\edge,\wavelength}^{1}(\xiparrand^1)}
\def \linksceMSIP {x_{\edge,\wavelength}^{t,(s,d)}(\xisupt)}
\def \wavelinksceMSIP {\mu_{\edge,\wavelength}^{t}(\xisupt)}
\def \wavelinkprevsceMSIP {\mu_{\edge,\wavelength}^{t-1}(\xisuptminusone)}

\def \numgrantedsceMSIP {z^{t,(s,d)}(\xisupt)}
\def \nodefunc {\delta}
\def \traffic {R}

\begin{document}
	
\RUNAUTHOR{Daryalal, Bodur, and Luedtke}

\RUNTITLE{Lagrangian Dual Decision Rules for MSMIP}

\TITLE{Lagrangian Dual Decision Rules for Multistage Stochastic Mixed-Integer Programming}

\ARTICLEAUTHORS{
	\AUTHOR{Maryam Daryalal}
	\AFF{Department of Decision Sciences, HEC Montréal, Montréal, Québec H3T 2A7, Canada, \EMAIL{maryam.daryalal@hec.ca}}
	\AUTHOR{Merve Bodur}
	\AFF{Department of Mechanical and Industrial Engineering, University of Toronto, Toronto, Ontario M5S 3G8, Canada, \EMAIL{bodur@mie.utoronto.ca}}
	\AUTHOR{James R. Luedtke}
	\AFF{Department of
		Industrial and Systems Engineering, University of Wisconsin, Madison, Wisconsin 53706, \EMAIL{jim.luedtke@wisc.edu}}
}
	
\ABSTRACT{
Multistage stochastic programs can be approximated by restricting policies to follow decision rules. Directly applying this idea to problems with integer decisions is difficult because of the need for decision rules that lead to integral decisions. In this work, we introduce \emph{Lagrangian dual decision rules} (LDDRs) for multistage stochastic mixed-integer programming ({\MSMIP}) which overcome this difficulty by applying decision rules in a Lagrangian dual of the \MSMIP. We propose two new bounding techniques based on stagewise (SW) and nonanticipative (NA) Lagrangian duals where the Lagrangian multiplier policies are restricted by LDDRs. We demonstrate how the solutions from these duals can be used to drive primal policies. Our proposal requires fewer assumptions than most existing {\MSMIP} methods. We compare the theoretical strength of the restricted duals and show that the restricted NA dual can provide relaxation bounds at least as good as the ones obtained by the restricted SW dual. In our numerical study \rev{on two problem classes, one traditional and one novel}, we observe that the proposed LDDR approaches yield significant optimality gap reductions compared to existing general-purpose bounding methods for {\MSMIP} problems.
}
	
	
\maketitle


\section{Introduction}\label{sec:intro} 
\begin{figure}[htb]
	\centering
	\includegraphics[scale=0.85]{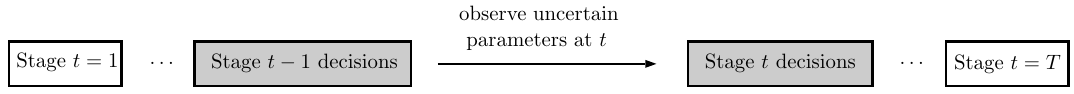}
	\caption{Dynamics of multistage stochastic programming}
	\label{fig:mssip-dynamics}
\end{figure}

\emph{Multistage stochastic mixed-integer programming} ({\MSMIP}) is a framework to model an optimization problem involving stochastic uncertainty, where the planning horizon is divided into multiple stages, decisions are made in each stage, and some of these decisions are constrained to be integer. The decisions in different stages cannot be made independently as they may impact subsequent stage decisions. The uncertainty is modeled as a stochastic process where the outcomes of random variables are observed over  stages.
In this setting, at each stage, the corresponding set of uncertain parameter values are observed, and based on this observation, the next stage decisions are made (see Figure \ref{fig:mssip-dynamics}). Therefore, as functions of random variables, decision variables at each  stage are random variables themselves.
However, they are \emph{nonanticipative}, i.e., they only depend on the history of observations, not future realizations. 
Thus,  in a multistage stochastic programming model, the solution is a \emph{policy} or \emph{decision rule} that maps all the past information to the current decisions to be made. 

{\MSMIP} naturally arises in many applications, such as unit commitment \citep{takriti1996stochastic}, capacity expansion \citep{rajagopalan1998capacity, ahmed2003multi, singh2009dantzig}, generation scheduling in hydro systems \citep{flatabo1998short, nowak2000stochastic, mo2001optimisation, flach2010long, helseth2015co}, batch sizing \citep{lulli2004branch}, airline fleet composition \citep{listes2005scenario}, transmission investment planning \citep{newham2007transmission}, transportation network protection  \citep{fan2010solving} and surgery planning \citep{gul2015progressive}.   

The majority of solution methods for {\MSMIP} require strong assumptions, including but not limited to stagewise independence, only right-hand-side uncertainty, binary state variables, and a finite (and not too large) scenario tree representation. Models lacking these conditions sometimes can be reformulated to satisfy the required assumptions, but at the expense of introducing (a potentially large number of) new variables and constraints. The most common assumption in the existing MSMIP methodologies is that the stochastic process is represented by a scenario tree.  
However, in general the size of a scenario tree required to obtain good quality solutions grows exponentially with the number of stages \citep{shapiro2005complexity}, which makes methods that rely on scenario tree models  unviable when the number of stages is beyond three or four. 
In this paper we propose new approximation approaches 
for {\MSMIP} problems that do not require these assumptions. These approaches are based on Lagrangian dual of an {\MSMIP} model and tractability is achieved by considering restricted forms of Lagrangian multipliers, i.e., forcing them to follow some \emph{decision rules}. 


\subsection{Related Literature}\label{sec:lit}

Even without integer decision variables, solving  an {\MSMIP} problem is theoretically and computationally challenging, due to high dimensional integration and the need to consider, while making current decisions, the (optimal) future decisions that will be made in response to the uncertain future trajectory of the stochastic process. Accordingly, all existing methods solve the problem by some form of approximation. There are three common \emph{approximation approaches} in the literature.

First is to model the underlying stochastic process in the form of a \emph{scenario tree} \citep{shapiro2009lectures}, which is the most common strategy in the existing solution methods for {\MSMIP}. Given a finite scenario tree approximation, the {\MSMIP} \rev{problem} is converted to a (very) large-scale, structured, deterministic problem. Techniques for solving such a scenario-tree based approximation include bounding techniques \citep{norkin1998branch,caroe1999dual,ahmed2003multi,alonso2003bfc,lulli2004branch,singh2009dantzig}, cutting-plane based methods \citep{guan2009cutting}, aggregation approaches \citep{sandikcci2014scalable}, nonanticipative Lagrangian dual approaches \citep{takriti1996stochastic,chen2002scenario}, and progressive hedging algorithms \citep{lokketangen1996progressive,listes2005scenario,fan2010solving,watson2011progressive,gul2015progressive,gade2016obtaining}. A limitation of the scenario tree approach is that, in order to obtain a good quality approximation of the stochastic process, the size of the tree in general needs to grow exponentially with the number of stages \citep{shapiro2005complexity}. 

Stochastic dual dynamic programming (SDDP), proposed by \cite{pereira1991multi}, is a leading solution method for solving multistage stochastic linear programming ({\MSLP}) problems as it is able to solve problems having an implicitly represented exponentially large scenario tree, under the assumption of stagewise independence. SDDP has been extended to {\MSMIP} by considering various approximations of the non-convex cost-to-go functions \citep{newham2007transmission,cen2012solving,cerisola2012stochastic, lohndorf2013optimizing, philpott2016midas}.
Proposing a new class of cutting planes, \cite{zou2017large} introduced stochastic dual dynamic integer programming (SDDiP) algorithm, an extension of SDDP which further assumes binary state decision variables. SDDiP is able to overcome some of its restrictions with various forms of reformulations, but at the expense of introducing new decision variables, thus increasing the size of the model. For instance, bounded integer state variables can be handled via binarization schemes. Most recently, \cite{ahmed2019stochastic} proposed stochastic Lipschitz dynamic programming for {\MSMIP} problems with general integer variables and Lipschitz cost-to-go functions, which uses Lipschitz cuts in its backward pass. The need to explicitly convexify the value function of many mixed-integer programming problems may potentially limit the scalability of this approach.

The second approximation approach, commonly known as the \emph{decision rule approach}, restricts the policies to follow a certain form, rather than restricting the form of the stochastic process. 
In the context of {\MSLP}, 
\cite{shapiro2005complexity} 
presented an upper bounding technique employing primal decision rules, where the decisions at each stage are restricted to be an affine function of observed random outcomes up to that stage. This yields a {static linear decision rule (LDR) policy}. On the other hand,  \cite{kuhn2011primal} provided a lower bounding technique  by applying LDRs to dual policies. \cite{bodur2018two} introduced a two-stage approximation by restricting only {state} variables (i.e., the ones  linking two consecutive stages together) to follow LDRs. By applying the so-called {two-stage LDRs} to primal and dual policies, the authors provided improved upper and lower bounds for {\MSLP} problems. 
Decision rules in the form of polynomial \citep{bampou2011scenario}, piecewise linear \citep{chen2008linear}, bilinear and trilinear \citep{georghiou2015generalized} functions have also been examined in the literature. 
In this \rev{approximation approach}, lower bounding techniques 
are limited to {\MSLP} problems. We
extend this line of research to {\MSMIP} by proposing the use of decision rules to obtain tractable approximations of {\it Lagrangian} duals of {\MSMIP} problems. 

The final approximation approach includes work that does not assume a scenario tree but relies instead on exploiting \emph{problem structure}.  \cite{brown2010information} study information relaxation for stochastic dynamic programs and, similar to our approach, penalize the violation of  {nonanticipativity} constraints (i.e., the ones that ensure consistency among the decisions when perfect information about the future is assumed). Their penalty is based on a different dual, however, and the approach requires determining problem-specific penalty functions that balance computational tractability with the strength of the obtained bound.  \cite{barty2010decomposition} propose dual approximate dynamic programming, which provides approximations for problems that can be decomposed into smaller, tractable problems when certain linking constraints are relaxed.


\subsection{Our Contributions}

We introduce \emph{Lagrangian dual decision rules} (LDDRs) to obtain bounds for {\MSMIP} problems. We design \emph{two new lower (upper) bounding techniques} for general {\MSMIP} problems with minimization (maximization) objective, i.e., with mixed-integer state and recourse variables, without restricting the form of the underlying stochastic process. These bounding techniques are based on two Lagrangian relaxations of the {\MSMIP} model: stagewise (SW) and nonanticipative (NA), where the state equations linking consecutive stages and the nonanticipativity constraints are relaxed, respectively. In order to obtain  tractable approximations, we restrict the associated Lagrangian multiplier policies to follow a decision rule determined by parameters which are optimized to achieve the best possible bound given the restriction. We compare the theoretical strength of the restricted duals and show that, when appropriately constructed, the restricted NA dual provides a relaxation bound at least as good as the bound obtained by the restricted SW dual.

We  also develop \emph{two LDDR driven} primal policies  which can be used to obtain upper (lower) bounds for problems with minimization (maximization) objective, that incorporate the Lagrangian multipliers obtained from the lower (upper) bounding methods. We perform an extensive computational study on  a stochastic multi-item lot-sizing problem with lag  where demands follow an autoregressive process, and find that our SW dual based primal method returns good quality solutions, while our NA dual driven primal method provides better solutions at the expense of extra computational effort. Putting the obtained lower and upper bounds together, we observe that the LDDR restricted NA dual approach yields significant optimality gap reductions compared to standard general purpose bounding methods for {\MSMIP} problems. \rev{Moreover, we study a novel multistage provisioning problem in optical networks, namely routing and wavelength assignment problem, and show the considerable merit of using our framework in reducing the optimality gap from classical bounding methods.}

Our approximation approach is general purpose and free of strong assumptions made in the literature such as stagewise independence or existence of a tractable-sized scenario tree representation. 
To the best of our knowledge, this is the first approach based on  decision rules that is capable of handling mixed-integer state variables. Our approach leads to subproblems in the form of deterministic mixed-integer programs, and thus can exploit the state of the art in efficiently solving them. Moreover, as the form of the restricted Lagrangian dual problems is a two-stage stochastic program, our approach enables application of theory and methods for solving two-stage stochastic integer \rev{programming models}.  

\rev{We summarize the main contributions of our work as follows:
\BI
\I \textbf{Restricted Lagrangian dual problems}. Our proposal to apply the decision rule restriction to the two considered Lagrangian dual problems enables the use of these dual problems to obtain bounds for {\MSMIP} problems without requiring a small-size scenario tree.  In contrast, solving the full SW or NA Lagrangian dual problems is only feasible for {\MSMIP} models in which the stochastic process is represented with a scenario tree where the total number of nodes is of manageable size, dramatically limiting both the number of stages and the number of branches per stage that can be used to approximate the underlying stochastic process. For  example, in a capacity planning problem  by \cite{singh2009dantzig} it was observed that the NA dual could not be solved in a reasonable time for an instance defined by a scenario tree with four stages and eight branches per stage. We remark that while dual decision rule restrictions have previously been applied for MSLPs \citep{kuhn2011primal,bodur2018two}, our work is the first to use such restrictions for MSMIP.
\I \textbf{Theoretical comparison of Lagrangian dual bounds}. Previous theory on the unrestricted Lagrangian SW and NA dual problems demonstrates that the NA dual is at least as strong as the SW dual. However, this theory does not apply when both dual problems are restricted as in our proposal. We extend this theory to demonstrate how the basis functions for the NA dual can be selected, relative to the ones for the SW dual, to assure that the LDDR-restricted NA dual bound is at least as good as that of the LDDR-restricted SW dual bound.
\I \textbf{Dual-driven primal policies}. We develop new primal policies that use information from the solution of the restricted Lagrangian dual problems to guide decision making at each decision stage.
\I \textbf{Empirical study}: We conduct computational experiments on two problem classes. We validate our theory that the bounds obtained by the restricted NA dual problem are better than those obtained by the restricted SW dual problem, and demonstrate that the former can improve significantly over the classical perfect information bound. We also find that the proposed dual-driven primal policies can outperform a policy driven based on solving deterministic approximations. Moreover, one of our applications contributes to the telecommunications literature by presenting the first {\MSMIP} model for a widely-studied provisioning problem in optical networks.
\EI
}

The remainder of this paper is organized as follows. In Section   \ref{sec:lagrange-duals} we formally state the {\MSMIP} problem and two  Lagrangian dual problems arising from it. In Section \ref{sec:restricted-duals} we introduce new relaxation/bounding methods for an {\MSMIP} model that are based on those Lagrangian relaxations. In Section \ref{sec:ub-policy} we provide feasible policies designed by using the information obtained from the restricted duals. In Section \ref{sec:results} we evaluate the proposed methodologies \rev{on two problem classes}.

{\bf Notation.} Random variables are represented with bold letters ($ \boldsymbol{\xi} $) while their observations are regular font $ (\xi) $. We use $ [a] := \{1,2, \dots, a\} $ and $ [a,b] := \{a, a+1, \dots, b\} $ for positive integers $a$ and $b$ (with $a\leq b $), and $ (.)^{\top}$ for the transpose operator.


\section{Problem Statement}\label{sec:lagrange-duals}

Let $ T $ denote the number of decision stages and 
$ \xisubtrand $ be a random vector at stage $ t \in [T]$ with outcomes $ \xisubtval \in \R^{\ell_t} $, and $ \xi_1=1$ (i.e., the first stage is deterministic).   The stochastic process is represented by $ \{\xiparrand_t\}_{t=1}^T $ having probability distribution $ \dist $ and support $ \Xi $.  By
$ \xisupt =(\xiparrand_1,\dots,\xiparrand_t)$, we denote the history of the process at stage $ t $. \rev{By including $\xi_1 = 1$ in the history, we make sure that the stage $t$ decisions are functions of the entire history, including the first deterministic stage.}

An {\MSMIP} problem can be modeled as follows
\bsub
\label{eqs:mssip}
\begin{alignat}{2}
	\min_{ \revtwo{\{x_t(\cdot) \}_{t \in \set{T}} }}  \ \ & \Exp_{\xiT}\left[\sum_{t\in \set{T}}c_t(\xisupt)^\top x_t(\xisupt)\right]  \label{eq:obj-stochastic-ip} \\
	\text{s.t.} \ \ & A_t(\xisupt)x_t(\xisupt)+B_t(\xisupt)x_{t-1}(\xisuptminusone) = b_t(\xisupt), && \quad  t\in \set{T}, \allxi \label{eq:const-state-equation} \\
	& x_t(\xisupt)\in X_t(\xisupt), && \quad  t\in \set{T}, \allxi \label{eq:const-recourse-constraint} 
\end{alignat}
\esub
where $ c_t  : \R^{\ell^t} \rightarrow \R^{\nt} $,
$ A_t  : \R^{\ell^t} \rightarrow \R^{\mstatet\times \nt} $,	
$ B_t  : \R^{\ell^t} \rightarrow \R^{\mstatet\times \ntminusone}$,	
$ b_t  : \R^{\ell^t} \rightarrow \R^{\mstatet} $ and $ \ell^t = \sum_{t'=1}^t\ell_{t'} $. 
$ x_t(\xisupt) \in \Z^{p_t} \times \R^{\nt - p_t} $ are the decision variables, i.e., nonanticipative policies, for $ t\in\set{T}$. Here and throughout the paper we adopt the convention $x_0(\xi^0) \equiv 0$ .   
$ \allxi $ means that the constraints are required to be almost surely satisfied with respect to $ \dist $.
Constraints \eqref{eq:const-state-equation} and \eqref{eq:const-recourse-constraint} are \emph{state} and \emph{recourse} constraints, respectively, where $ X_t(\xisupt):=\{x \in \Z^{p_t} \times \R^{\nt - p_t} : C_t(\xisupt)x\geq d_t(\xisupt)\}$ with 	
$ C_t  : \R^{\ell^t} \rightarrow \R^{\mrecourset\times \nt} $ and	
$ d_t : \R^{\ell^t} \rightarrow \R^{\mrecourset}$.
\rev{

Throughout this work, we make the following assumptions: 
\begin{assumption}\label{ass:rel-recourse} Model \eqref{eqs:mssip} has relatively complete recourse, i.e., \revtwo{ for each stage $t \in [1,T]$, and every possible history $\xisupt$ and sequence of feasible decisions $x_{s}(\mathbf{\xi}^s)$ at $s\in[1,t-1]$, there always exists a feasible decision $x_t(\xisupt)$.}
\end{assumption}
\begin{assumption}\label{ass:bounded}
\revtwo{The set $X_t(\xisupt)$ is compact for all $ t\in\set{T} $, $ \allxi $ and $\Exp[\dia{X_t(\xisupt)}]$ is finite, where $\dia{X_t(\xisupt)}:= \max\{ \| x - y\|_2 : x,y \in X_t(\xisupt)\}$ is the diameter of $X_t(\xisupt)$.}
\end{assumption}
Assumption \ref{ass:rel-recourse} is common in the stochastic programming literature. It is naturally satisfied in many applications, even in the stronger form of complete recourse \citep{rockafellar1976stochastic, LULLI2006879}, e.g., \cite{zou2017stochastic} and \cite{ahmed2019stochastic} both assume complete continuous recourse. 
Otherwise, it can be achieved by modifying the model to include slack variables that measure constraint violations \citep{shapiro2005complexity, zou2017stochastic} which are then penalized in the objective. Assumption \ref{ass:bounded} refers to the boundedness of the feasibility set, which is often the case in practice. If it is not readily implied by the problem definition, it can be enforced by adding sufficiently large (small) upper (lower) bounds on decision variables.
}

We next describe two Lagrangian duals for {\MSMIP}, namely SW dual \citep{rosa1996augmented} and NA dual \citep{rockafellar1991scenarios}.


\subsection{Stagewise Lagrangian Dual}
In the SW dual, the constraints linking consecutive stages are relaxed, namely the \emph{state equations}. Let $ \pi_t(\xisupt) \in \R^{\mstatet}, t\in\set{T} $ be the dual variables associated with constraints \eqref{eq:const-state-equation}. For fixed dual functions $ \pi_t: \R^{\ell^t} \rightarrow \R^{\mstatet}, t\in\set{T} $, the SW Lagrangian relaxation problem is defined as
\bsub
\label{eqs:sw}
\begin{align}
\hspace{-5mm}\mathfrak{L}^{\text{SW}}(\pi_1,\dots,\pi_T) = \nonumber
\\
\min_{ \revtwo{\{x_t(\cdot) \}_{t \in \set{T}}} }  \ & \Exp_{\xiT} \left[\sum_{t\in \set{T}} c_t(\xisupt)^\top x_t(\xisupt) + \pi_t(\xisupt)^\top \left(A_t(\xisupt)x_t(\xisupt)+B_t(\xisupt)x_{t-1}(\xisuptminusone) - b_t(\xisupt)\right)  \right]  \label{eq:obj-lagrange-sw-t} \\
\text{s.t.} \  & x_t(\xisupt)\in X_t(\xisupt), \quad  t\in \set{T}, \allxi, \label{eq:const-lagrange-sw-recourse-cons}
\end{align}
\esub
which decomposes by \revtwo{ stage $ t $ and sample path $ \xi^T $, thus can be solved by solving an independent deterministic mixed-integer program (MIP) per stage $t$ (stage-decomposable) and scenario $\xi^T$ (path-decomposable)}. As the Lagrangian relaxation problem provides a valid lower bound on the optimal value of the {\MSMIP} problem, the stagewise Lagrangian dual problem aims to find the dual functions providing the best bound
\begin{equation}\label{eq:sw-dual}
\nu^{\text{SW}} := \max_{ \{\pi_t(\cdot) \}_{t \in \set{T}} } \mathfrak{L}^{\text{SW}}(\pi_1,\dots,\pi_T).
\end{equation}


\subsection{Nonanticipative Lagrangian Dual}
The NA dual is based on a reformulation of the {\MSMIP} problem where we create a copy of every decision variable for every realization and explicitly enforce nonanticipativity. We introduce the copy variables $ y(\xiT)=(y_1(\xiT),\dots,y_T(\xiT)) $ as perfect information variables, meaning that they depend on the entire sample path $\xiT= (\xiparrand_1,\dots,\xiparrand_T) $. For every sample path $\xiT \in \Xi$, we define the set
\[ Y(\xiT) = \bigl\{y \in \Z^{p_t} \times \R^{\nt - p_t}: A_t(\xisupt)y_t+B_t(\xisupt)y_{t-1} = b_t(\xisupt), y_t\in X_t(\xisupt), t\in\set{T} \bigr\}
. \]
Then, the {\MSMIP} problem \eqref{eqs:mssip} can be reformulated as
\bsub
\begin{alignat}{2}
\min_{ \revtwo{\{y_t(\cdot) \}_{t \in \set{T}}} }  \ \ & \Exp_{\xiT}\left[\sum_{t\in \set{T}}c_t(\xisupt)^\top y_t(\xiT)\right]  \label{eq:obj-reform-stochastic-ip} \\
\text{s.t.} \ \ & y(\xiT) \in Y(\xiT), && \quad  \allxi \label{eq:const-perfect-info}\\
& y_t(\xiT) = \Exp_{\xiprimeT}\left[y_t(\xiprimeT) \big| \xiprimesupt = \xisupt \right], && \quad t\in \set{T}, \allxi. \label{eq:const-reform-nonanticipative} 
\end{alignat}
\esub
Constraints  \eqref{eq:const-reform-nonanticipative} are \emph{nonanticipativity constraints} which make sure that for every partial realization of a sample path $\xisupt$ at stage $t$, the decisions made at stage $t$ are consistent (i.e., the decisions made in all sample paths $\xiT$ that share the history $\xisupt$ are the same).
Associating the dual functions $  \gamma_t(\cdot) \in \Gamma_t, t\in \set{T} $, where $ \Gamma_t =\{\gamma_t : \R^{\ell^T} \rightarrow \R^{\nt}\ |\ \Exp[\gamma_{t}(\xiT)] < \infty \}$, the nonanticipativity constraints are relaxed to obtain the NA Lagrangian dual problem
\bsub
\label{eqs:na}
\begin{align}
\mathfrak{L}^{\text{NA}}(\gamma_1,\dots,\gamma_T) = \nonumber\\
\min_{ \revtwo{\{y_t(\cdot) \}_{t \in \set{T}}} } \ \ & \Exp_{\xiT}\left[\sum_{t\in \set{T}}c_t(\xisupt)^\top y_t(\xiT) + \gamma_t(\xiT)^\top \left(y_t(\xiT) - \Exp_{\xiprimeT}[y_t(\xiprimeT) | \xiprimesupt = \xisupt]\right)\right]  \label{eq:obj-lagrange-na-t} \\
\text{s.t.} \ \ & y(\xiT) \in Y(\xiT), \quad \allxi. \label{eq:const-lagrange-na-perfect-info}
\end{align}
\esub
The objective function in \eqref{eq:obj-lagrange-na-t} can be simplified using the following lemma whose proof can be found in 
the Appendix. \rev{This lemma allows us to swap the Lagrangian multipliers $\gamma_t(\xiT)$ and decision variables inside the expectation $y_t(\xiprimeT)$, which makes it possible to compute the coefficients ahead of the time.}
 \begin{lemma}\label{lem:na-reform}
 	Assume that \revtwo{$\Exp[\dia{X_t(\xisupt)}]$} is finite. Then, the following equality holds
 	\begin{align}
		& \Exp_{\xiT}\left[\sum_{t\in \set{T}}c_t(\xisupt)^\top y_t(\xiT) + \gamma_t(\xiT)^\top \left(y_t(\xiT) - \Exp_{\xiprimeT}[y_t(\xiprimeT) | \xiprimesupt  = \xisupt]\right)\right]  =  \label{eqq:lemma1} & \\
		& \Exp_{\xiT}\left[\sum_{t\in \set{T}}\left(c_t(\xisupt)  + \gamma_t(\xiT) - \Exp_{\xiprimeT}[\gamma_t(\xiprimeT) | \xiprimesupt = \xisupt]\right)^\top y_t(\xiT) \right] .\nonumber&
 	\end{align}
 \end{lemma}
\noindent
Using Lemma \ref{lem:na-reform}, the nonanticipative Lagrangian relaxation for  fixed dual functions $ \gamma_t : \R^{\ell^T} \rightarrow \R^{\nt}, t\in\set{T} $ can be written as
\begin{align*}
\mathfrak{L}^{\text{NA}}(\gamma_1,\dots,\gamma_T) = \min_{ \revtwo{\{y_t(\cdot) \}_{t \in \set{T}}} } \ \ & \Exp_{\xiT}\left[\sum_{t\in \set{T}}\left(c_t(\xisupt)  + \gamma_t(\xiT) - \Exp_{\xiprimeT}[\gamma_t(\xiprimeT) | \xiprimesupt = \xisupt]\right)^\top y_t(\xiT) \right] 
 \\
\text{s.t.} \ \ & y(\xiT) \in Y(\xiT), \quad \allxi, 
\end{align*}
\noindent
which is decomposable by sample path, but not by stage. \revtwo{Specifically,  $\mathfrak{L}^{\text{NA}}(\gamma_1,\dots,\gamma_T)$ can be evaluated by solving a separate $T$-period deterministic mixed-integer program for each sample path $\xiT$.}
Finally the NA Lagrangian dual problem is
\begin{equation}\label{eq:na-dual}
\nu^{\text{NA}} := \max_{\{ \gamma_t(\cdot)\in \Gamma_t \}_{t \in [T]}} 
\mathfrak{L}^{\text{NA}}(\gamma_1,\dots,\gamma_T).
\end{equation}


\subsection{Primal Characterizations and Bound Comparison}\label{sec:bound-comp-regular}
Using Lagrangian duality theory,  
 primal characterizations of the Lagrangian duals can be obtained, which  in turn can be used to compare their strength \citep{dentcheva2004duality}.
For the SW Lagrangian dual problem, the primal characterization is as follows
\bsub
\label{eq:sw-primal-char}
\begin{alignat}{2}
\nu^{\text{SW}} = \min_{ \revtwo{\{x_t(\cdot) \}_{t \in \set{T}} }} \ \ & \Exp_{\xiT}\left[\sum_{t\in \set{T}}c_t(\xisupt)^\top x_t(\xisupt)\right]  \\
\text{s.t.} \ \ & A_t(\xisupt)x_t(\xisupt)+B_t(\xisupt)x_{t-1}(\xisuptminusone) = b_t(\xisupt), && \quad t\in \set{T}, \allxi  \\
& x_t(\xisupt)\in \conv{X_t(\xisupt)}, && \quad  t\in \set{T}, \allxi. 
\end{alignat}
\esub
That is, for each stage and sample path, the feasible set of the recourse problem on that stage is relaxed and replaced with its convex hull. 
The primal characterization of the NA Lagrangian dual problem is given below
\begin{alignat*}{2}
\nu^{\text{NA}} = \min_{ \revtwo{\{y_t(\cdot) \}_{t \in \set{T}} }} \ \ & \Exp_{\xiT}\left[\sum_{t\in \set{T}}c_t(\xisupt)^\top y_t(\xiT)\right]   \\
\text{s.t.} \ \ & y_t(\xiT)\in \conv{Y(\xiT)},  && \quad \allxi\\
& y_t(\xiT) = \Exp_{\xiprimeT}[y_t(\xiprimeT) | \xiprimesupt = \xisupt], && \quad t\in \set{T}, \allxi. 
\end{alignat*}
That is, for every sample path, the feasible set of the $T$-stage deterministic problems is relaxed and replaced with its convex hull.

Using these characterizations, \cite{dentcheva2004duality} show, in the case when the stochastic process is represented by a finite scenario tree, that the NA dual is not worse than the SW dual, i.e.,
\begin{equation*} 
\nu^{\text{SW}} \leq \nu^{\text{NA}}.
\end{equation*}


\section{Lagrangian Dual Decision Rules for {\MSMIP}}\label{sec:restricted-duals}
The SW dual functions $ \pi_t  $ and the NA dual functions $ \gamma_t $  
are policies which map every possible history of observations to a dual decision vector, making direct solution of the
respective Lagrangian dual problems \eqref{eq:sw-dual} and \eqref{eq:na-dual} intractable in general. We propose 
restricting these dual multipliers to follow decision rules, referred to as LDDRs. As such, we  obtain restricted problems that have finitely many decision variables. In what follows, we first explain the LDDR approach for the SW and NA duals (Sections  \ref{sec:restricted-sw} and \ref{sec:restricted-na}), then provide their primal characterizations for a strength comparison as in the unrestricted case (Section \ref{sec:restricted-primalchar}), and lastly provide an algorithmic framework for their solutions (Section \ref{sec:SAA}). 


\subsection{Restricted Stagewise Lagrangian Dual}\label{sec:restricted-sw}
For $t \in [T]$, we restrict the dual variables $ \pi_t(\xisupt) $ to follow an LDDR by enforcing 
$$\pi_t(\xisupt) =  \Phi_{t}(\xisupt) \ldrvec_t$$
where $ \Phi_t: \R^{\ell^t} \rightarrow \R^{\mstatet \times \Kt} $ are the set of \emph{basis functions} and $ \ldrvec_t \in \R^{\Kt}$ is a vector of \emph{LDDR decision variables} (i.e., the weights associated with the basis functions). For the ease of presentation, we use the matrix form of basis function outputs -- i.e., one can think of $\Phi_t$ as consisting of $\Kt$ basis functions, each of which maps the history $\xisupt$ to a vector of size $\mstatet$, the number of state equations in stage $t$. The set of basis functions is a model choice, and so $ \Kt $ depends on this choice.  We define the \emph{LDDR-restricted SW Lagrangian dual} problem which aims to find the optimal choice of LDDR variable values to maximize the obtained lower bound
\begin{equation}\label{eq: ldr-lagrange-dual}
\nu_R^{\text{SW}} := \max_{ \{\ldrvec_t\}_{t \in \set{T}}} \mathfrak{L}^{\text{SW}}(\Phi_1\ldrvec_1,\dots,\Phi_T\ldrvec_T) =
 \max_{\{\ldrvec_t\}_{t \in \set{T}}} \Exp_{\xiT}\left[\sum_{t\in\set{T}}\mathcal{L}^{\text{SW}}_t(\ldrvec_{t},\xisupt) - \ldrvec^\top_{t} \Phi^\top_{t}(\xisupt)\ b_t(\xisupt)\right]
\end{equation}
where  
$\displaystyle \mathcal{L}^{\text{SW}}_t(\ldrvec,\xisupt) := \min_{ \revtwo{\{x_t \}_{t \in \set{T}} }} \bigg\{ \Big(c_t(\xisupt) + \ldrvec^\top_{t} \Phi^\top_{t}(\xisupt)\ A_t(\xisupt) +\ldrvec^\top_{t+1}  \Exp_{\xiprimeT}\big[\ \Phi^\top_{t+1}(\xiprimesuptplusone)\  B_{t+1}(\xiprimesuptplusone) \ |  \ \xisupt \big]\Big)^\top x_t: x_t\in X_t(\xisupt) \bigg\}. $ \rev{Observe that, for fixed $\xisupt$, calculating the coefficient of $ \ldrvec_{t+1} $ requires evaluating the conditional expectation $\Exp_{\xiprimeT}[\Phi^\top_{t+1}(\xiprimesuptplusone)B_{t+1}(\xiprimesuptplusone) | \xisupt] $. The tractability of this calculation depends on the form of $\Phi_t$ and the dependence structure of the stochastic process. For instance, they can be directly calculated if $\Phi_t$ consists of affine functions and the conditional distribution of $\xiprimesuptplusone$ given $\xisupt$ is known, such as when the random variables follow an autoregressive process. In general, these coefficients can be estimated via sampling.}


\subsection{Restricted Nonanticipative Lagrangian Dual}\label{sec:restricted-na}
Letting $ \Psi_t: \R^{\ell^T} \rightarrow \R^{\nt \times \Kprimet} $ for $t \in [T]$ to be a set of basis functions, we restrict the dual variables $ \gamma_t(\xiT) $ to follow an LDDR
$$ \gamma_t(\xiT) = \ \Psi_t(\xiT)\naldrvec_t,$$
where  $ \naldrvec_t\in \R^{\Kprimet}, t\in \set{T}, $ is the vector of LDDR decision variables. Then, we obtain the \emph{LDDR-restricted NA Lagrangian dual} problem as
\begin{equation}\label{eq: ldr-na-lagrange-dual}
\nu_R^{\text{NA}} := \max_{ \{\alpha_t\}_{t \in \set{T}} } \mathfrak{L}^{\text{NA}}( \Psi_1\naldrvec_1,\dots, \Psi_T\naldrvec_T) =
 \max_{\{\alpha_t\}_{t \in \set{T}}} \Exp_{\xiT}\left[ \mathcal{L}^{\text{NA}}(\naldrvec,\xiT) \right],
\end{equation}
where  
\begin{align}
\displaystyle \mathcal{L}^{\text{NA}}(\naldrvec,\xiT) = & \min_{ \revtwo{\{y_t \}_{t \in \set{T}} }} \left\lbrace \sum_{t\in \set{T}} \left(c_t(\xisupt)  + \Psi_t(\xiT)\naldrvec_t - \Exp_{\xiprimeT}[\Psi_t(\xiprimeT)\naldrvec_t | \xiprimesupt = \xisupt]\right)^\top y_t : y \in Y(\xiT) \right\rbrace \nonumber\\
= & \min_{ \revtwo{\{y_t \}_{t \in \set{T}} }} \left\lbrace \sum_{t\in \set{T}} \left(c_t(\xisupt)  + \left(\Psi_t(\xiT) - \Exp_{\xiprimeT}[\Psi_t(\xiprimeT) | \xiprimesupt = \xisupt]\right) \naldrvec_t\right)^\top y_t : y \in Y(\xiT) \right\rbrace.\label{eq:na_subpproblem}
\end{align}
\rev{Since basis functions $\Psi_t(\xiT)$ are functions of the entire sample path, depending on their structure we might need to evaluate the $h$-step ahead conditional expectations, for $1\leq h \leq T-t$, while for the SW restricted dual we have $h=1$. For example, if $\Psi_{t\kprimeindex}(\xiT) = \xiparrand^{t+h}$, then $\Exp_{\xiprimeT}[\xiparrand^{t+h} | \xiprimesupt = \xisupt]$ needs to either be calculated by an analytical form (if possible), or estimated by means of sampling. See Sections \ref{sec:data-instances} and \ref{sec:data-bf-RWA} for some examples of computing the explicit form of the above conditional expectations.}


\subsection{Primal Characterizations and Bound Comparison}
\label{sec:restricted-primalchar}

In this section we assume that the stochastic process is represented by a finite scenario tree, implying that for every stage $t \in [T]$, the random vector $\xiparrand_t$ has a discrete distribution with finite support, and also that the set of possible sample paths, $\Xi$, is finite. We make this assumption in this section only for the sake of simplicity in deriving primal characterizations for the proposed restricted Lagrangian duals, and note that we do not require the size of the scenario tree representation to be tractable as this assumption is used only for theoretical \rev{analysis}, not for computational purposes. We emphasize that our overall methodological framework does not restrict the stochastic process to a finite scenario tree model, and rather relies on sampling to approximately solve the restricted approximations. The finite scenario tree assumption of this section is in line with the comparison of the bounds from unrestricted duals by \cite{dentcheva2004duality}.  

We next derive a general primal characterization of the Lagrangian dual of a MIP when dual multipliers are restricted to a linear form. Since an {\MSMIP} is a large-scale (structured) MIP under the finite scenario tree assumption, this analysis is sufficient to obtain the primal characterizations of our restricted duals of the {\MSMIP} problem. 

Consider a MIP $ \min \{c^\top x : Dx= d, \ x\in X\} $ where $X \subseteq \mathbb{R}^n$ is a set defined by linear constraints and integer constraints on some of the decision variables, $D \in \mathbb{R}^{m \times n}$, $c \in \mathbb{R}^n$, and $d \in \mathbb{R}^m$. 
For $\lambda \in \mathbb{R}^m$ define
\begin{equation}
\label{eq:zld} 
z(\lambda) = \min_{ \revtwo{x}} \{ c^\top x + \lambda^\top (d-Dx) : x \in X \}. 
\end{equation}
The standard Lagrangian dual is the problem $z^{\text{LD}} = \displaystyle \max_{\lambda} z(\lambda)$, and its primal characterization is
\begin{equation}
\label{eq:ldprim}
 z^{\text{LD}} = \min_{ \revtwo{x}} \{ c^\top x : x \in \conv{X}, Dx = d \} .
\end{equation}
For a given matrix $G \in \mathbb{R}^{m \times K}$, define the the {\it restricted} Lagrangian dual
 \begin{equation}
 \label{eq:rld}
 z^{\text{RLD}} \ = \ \max_{\lambda, \alpha} \{ z(\lambda) : \lambda = G\alpha \} . 
\end{equation}
The following Lemma provides a primal characterization of this restricted dual problem.

 \begin{lemma}\label{lem:primal-char}
 The restricted Lagrangian dual satisfies
 	\begin{alignat*}{2}
 	z^{\text{RLD}} \ =   \min_x & \ c^\top x  \\
    \ \text{s.t.} \ &  x\in\conv{X} \\
 	 	&  G^\top  (Dx - d) = 0.
 	\end{alignat*}
 \end{lemma}
 In other words, in comparison to the unrestricted Lagrangian dual \eqref{eq:ldprim}, the restricted Lagrangian dual still replaces the set $X$ with its convex hull, but the constraints $Dx = d$ are relaxed to $G^{\top} (Dx - d) = 0$. The proof  of Lemma \ref{lem:primal-char} can be found in
 the Appendix. 
 Using this lemma, we derive the primal characterizations of the restricted Lagrangian duals of the {\MSMIP}.

As the stochastic process is assumed to be represented by a finite scenario tree, the probability of each sample path can be assumed to be positive. In order to apply Lemma \ref{lem:primal-char}, we scale the state equations in the {\MSMIP} formulation \eqref{eqs:mssip} by the corresponding probabilities of these scenarios denoted by $ p(\xiTval) $
\bsub
\begin{alignat*}{2}
	\min_{ \revtwo{\{x_t(\cdot) \}_{t \in \set{T}} }}  \ \ & \Exp_{\xiT}\left[\sum_{t\in \set{T}}c_t(\xisupt)^\top x_t(\xisupt)\right]  \\
	\text{s.t.} \ \ & p(\xiTval) (A_t(\xisuptval)x_t(\xisuptval)+B_t(\xisuptval)x_{t-1}(\xisuptminusoneval)) = p(\xiTval)b_t(\xisuptval), && \quad  t\in \set{T}, 	\allxiscen
	 \\
	& x_t(\xisuptval)\in X_t(\xisuptval), && \quad  t\in \set{T}, \allxiscen.
\end{alignat*}
\esub
Recall the LDDR restrictions for the duals
$$ \pi_t(\xisuptval) -  \Phi_t(\xisuptval)\beta_t = 0, \quad  t\in \set{T}, \allxiscen. $$
These are the new set of constraints added to the SW Lagrangian dual problem \eqref{eq:sw-dual}. Consequently, using Lemma \ref{lem:primal-char}, the primal characterization of the restricted SW Lagrangian dual problem is obtained as
\bsub
\label{eqs:SWprimalchar}
\begin{alignat}{2}
\hspace*{-0.19cm} \nu_R^{\text{SW}} = \min_{ \revtwo{\{x_t(\cdot) \}_{t \in \set{T}} }}  & \Exp_{\xiT}\left[\sum_{t\in \set{T}}c_t(\xisupt)^\top x_t(\xisupt)\right]   \\
\text{s.t.}\ \ & x_t(\xisuptval)\in \conv{X_t(\xisuptval)}, &\quad&  t\in \set{T}, \allxiscen \label{eq:const-strength-ldr-recourse} \\
& \Exp_{\xiT}\left[  \Phi_{t}(\xisupt)^\top \left(A_t(\xisupt)x_t(\xisupt)+B_t(\xisupt)x_{t-1}(\xisuptminusone) - b_t(\xisupt)\right)\right] = 0, &&   t\in \set{T}.  \label{eq:const-strength-ldr-state}
\end{alignat}
\esub 
Note that, for each $ t\in\set{T} $, there are exactly $\Kt$ expected value constraints, one for each basis function vector in the matrix $\Phi_t$.
Compared to the primal characterization of the unrestricted SW dual \eqref{eq:sw-primal-char}, while the recourse problem feasible sets are still convexified (as in \eqref{eq:const-strength-ldr-recourse}),  
the state equations are not almost surely satisfied, 
rather for every basis function, their expectation is enforced to be equal to zero. Therefore, unlike the unrestricted version \eqref{eq:sw-primal-char}, even if 
the original problem is convex (such as an {\MSLP}), a duality gap might exist.


Using the same approach, we obtain the following primal characterization for the restricted NA Lagrangian dual problem
\bsub
\label{eqs:NAprimalchar}
\begin{alignat}{2}
\nu_R^{\text{NA}} = \min_{ \revtwo{\{y_t(\cdot) \}_{t \in \set{T}} }} \ \ & \Exp_{\xiT}\left[\sum_{t\in \set{T}}c_t(\xisupt)^\top y_t(\xiT)\right]   \\
\text{s.t.} \ \ & y_t(\xiTval)\in \conv{Y(\xiTval)}, &\quad& \allxiscen \label{eq:const-2-na-ldr-primal} \\
& \Exp_{\xiT}\Big[ \Psi_{t}(\xiT)^\top \big(y_t(\xiT) - \Exp_{\xiprimeT}[y_t(\xiprimeT) | \xiprimesupt = \xisupt]\big) \Big] = 0, && t\in \set{T}.  \label{eq:const-3-na-ldr-primal}
\end{alignat}
\esub
Here, constraint sets $ Y(\xiTval) $ are replaced with their convex hull, but the NA constraints are enforced in expectation rather than for almost every $ \xiT$.

As mentioned in Section \ref{sec:bound-comp-regular}, for the unrestricted duals, it is proven that 
$ \nu^{\text{SW}} \leq \nu^{\text{NA}} $. The same result does not immediately hold for the restricted duals due to the
relaxed constraints. However, we prove that when the basis functions for the restricted NA dual are \rev{properly} selected, a similar inequality holds for the restricted duals.
Let $ \Phi_t=(\Phi_{t1},\dots,\Phi_{t \Kt}) $ and $ \Psi_t=(\Psi_{t1},\dots,\Psi_{t \Kprimet}) $ be the basis functions 
used in the restricted SW and NA duals, respectively. Note that 
$ \Phi_t: \R^{\ell^t} \rightarrow \R^{\mstatet \times \Kt} $ and $ \Psi_t: \R^{\ell^T} \rightarrow \R^{\nt \times \Kprimet}$  are functions mapping a set of random
variables to a column vector for each basis function indexed by $k$ \rev{and $\kprimeindex$, respectively}.
\rev{
\begin{theorem}\label{thm:sw-na-compare}
    Given a set of basis functions $ \Phi_t=(\Phi_{t1},\dots,\Phi_{t \Kt}), t \in [T] $ used in restricted SW dual, there exists $ \Psi_t=(\Psi_{t1},\dots,\Psi_{t \Kprimet}), t\in[T] $ for the restricted NA dual such that:
    \begin{equation*} 
	\nu_R^{\text{SW}} \leq \nu_R^{\text{NA}}.
	\end{equation*}
\end{theorem}
}

\proof{Proof}
\rev{
In the following constructive proof, we first build a set of basis functions $\Psi_t$ given $\Phi_t$. Then we show that having such $\Psi_t$ guarantees $\nu_R^{\text{SW}} \leq \nu_R^{\text{NA}}$.

For each $ t\in \set{T} $, we define the set of basis functions $\Psi_t$, with cardinality $\Ksetprime_t = \Kpar_t + \Kpar_{t+1}$ for $t<T$ and $\Ksetprime_T = \Kpar_T$. For $ k\in\set{\Kt} $:
	\begin{enumerate}
		\item $ \Psi_{t \kprimeindex}(\xiT)^\top:=\Phi_{tk}(\xisupt)^\top A_t(\xisupt) $ for every $\xiT \in \Xi$, where $\kprimeindex=k$,
		\item If $t > 1$, then $\Psi_{{t-1} \kprimeindex}(\xiT)^\top :=\Phi_{tk}(\xisupt)^\top B_t(\xisupt) $ for every $\xiT \in \Xi$, where $\kprimeindex = \Kt + k$.
	\end{enumerate}
Using these basis functions $\Psi_t$,} let $ y^\star(\xiT) $ be an optimal solution of the restricted NA dual \eqref{eqs:NAprimalchar} with the optimal value of $ \nu^{\text{NA}}_R$. The existence of $ y^\star(\xiT) $ is guaranteed due to our initial assumption that the stochastic process is modeled as a finite scenario tree, as well as the boundedness of the feasible region.  We construct a feasible solution to the restricted SW dual \eqref{eqs:SWprimalchar} with the same objective value $ \nu^{\text{NA}}_R $, which demonstrates the desired inequality.
We claim that $ \hat{x}_t(\xisupt) =  \Exp_{\xiprimeT}[y^{\star}_t({\xiprimeT}) | \xiprimesupt = \xisupt]$ satisfies these conditions.

We  first show that at $ \hat{x}_t(\xisupt) $ the objective function value of the restricted SW dual \eqref{eqs:SWprimalchar} is equal to $ \nu^{\text{NA}}_R $. Indeed, the objective evaluates to 
\bsub
\begin{alignat}{4}
&\Exp_{\xiT}\left[\sum_{t\in \set{T}}c_t(\xisupt)^\top \hat{x}_t(\xisupt) \right]  & = &  \Exp_{\xiT}\left[\sum_{t\in \set{T}}c_t(\xisupt)^\top \Exp_{\xiprimeT}[y^{\star}_t(\xiprimeT) | \xiprimesupt = \xisupt]\right]  &  &  \nonumber \\\
&& =&  \sum_{t\in \set{T}} \Exp_{\xisupt}\left[\Exp_{\xiprimeT}[c_t(\xiprimesupt)^\top y^{\star}_t(\xiprimeT) | \xiprimesupt = \xisupt]\right] &  & \qquad{\color{gray} } \label{eq:obj-equality-2}\\
&& = & \sum_{t\in \set{T}} \Exp_{\xiprimeT}[c_t(\xiprimesupt)^\top y^{\star}_t(\xiprimeT) ]  & = & \qquad \nu^{\text{NA}}_R \label{eq:obj-equality-3}
\end{alignat}
\esub
where \eqref{eq:obj-equality-2} follows because for a fixed $ \xiT $, $c_t(\xisupt) $ is a constant that can be brought into the inside conditional expectation, and by swapping the order of expectation, the first equality in \eqref{eq:obj-equality-3} follows from the fact that, for random variables $Z$ and $Y$,  $\Exp_Y[\Exp_Z[Z|Y]]= \Exp[Z]$, and the second equality in \eqref{eq:obj-equality-3} follows because $\xiT$ and $\xiprimeT$ have the same distribution.

We next verify that $ \hat{x}_t(\xisupt) $ satisfies constraint \eqref{eq:const-strength-ldr-recourse}. By \eqref{eq:const-2-na-ldr-primal} and the definition of $Y(\xiprimeT)$, $ y^{\star}_t(\xiprimeT) \in X_t(\xisupt) $ for each $t \in T$ and $\xiprimeT \in \Xi$ such that $\xiprimesupt=\xisupt$. Since $\hat{x}_t(\xisupt) =  \Exp_{\xiprimeT}[y^{\star}_t(\xiprimeT) | \xiprimesupt = \xisupt] $ is a convex combination of such $ y_t^{\star}(\xiprimeT)$ variables, it follows that $\hat{x}(\xisupt) \in \conv{X_t(\xisupt)} $.

Finally, we verify that $ \hat{x}_t(\xisupt) $ satisfies constraint  \eqref{eq:const-strength-ldr-state}.
	We evaluate the left-hand side of \eqref{eq:const-strength-ldr-state} for a fixed $t \in [T]$ at the defined solution $\hat{x}$. First break the expectation into three terms
	\begin{align}
	&\Exp_{\xiT}\left[  \Phi_{tk}(\xisupt)^\top
	\left(A_t(\xisupt)\hat{x}_t(\xisupt)+B_t(\xisupt)\hat{x}_{t-1}(\xisuptminusone) - b_t(\xisupt)\right)\right] =
	\label{eq:lhs} \\
	&\Exp_{\xiT}\left[\Phi_{tk}(\xisupt)^\top A_t(\xisupt)\hat{x}_t(\xisupt) \right] +
	\Exp_{\xiT}\left[\Phi_{tk}(\xisupt)^\top B_t(\xisupt)\hat{x}_{t-1}(\xisuptminusone) \right] -
	\Exp_{\xiT}\left[\Phi_{tk}^\top(\xisupt) b_t(\xisupt)\right]. \label{eq:break}
	\end{align}
	 Substituting $ \hat{x}_t(\xisupt) =
	\Exp_{\xiprimeT}[y^{\star}_t(\xiprimeT) |\xiprimesupt = \xisupt]$ in the first term of \eqref{eq:break}, yields
	$$
	\Exp_{\xiT}\left[\Phi_{tk}(\xisupt)^\top A_t(\xisupt)\Exp_{\xiprimeT}[y^{\star}_t(\xiprimeT) | \xiprimesupt = \xisupt] \right]. 
	$$
	 \rev{By construction,} there exists some $ \kprimeindex\in\set{\Kprimet} $ with $\Psi_{t \kprimeindex}(\xiT)^\top=\Phi_{tk}(\xisupt)^\top A_t(\xisupt)$, 
	thus constraint \eqref{eq:const-3-na-ldr-primal} leads to the following equality
	$$\Exp_{\xiT}\left[\Phi_{tk}(\xisupt)^\top A_t(\xisupt)\Exp_{\xiprimeT}[y^{\star}_t(\xiprimeT) | \xiprimesupt = \xisupt] \right]  = \Exp_{\xiT}\left[\Phi_{tk}(\xisupt)^\top A_t(\xisupt)y^{\star}_t(\xiT) \right].$$
	 This follows by substituting the vector $ \Phi_{tk}(\xisupt)^\top A_t(\xisupt) $ on the left-hand-side expression  
	 with its equivalent vector $ \Psi_{t \kprimeindex}(\xiT)^\top $ and using \eqref{eq:const-3-na-ldr-primal} to show the equality
	 of $ \Exp_{\xiT}[\Psi_{t \kprimeindex}(\xiT)^\top \Exp_{\xiprimeT}[{y^*_t(\xiprimeT)_j}|\xiprimesupt =\xisupt]] =
	 \Exp_{\xiT}[\Psi_{t \kprimeindex}(\xiT)^\top(y^*_t(\xiT))_j]$. By applying the same argument to the second term of \eqref{eq:break} while employing constraint \eqref{eq:const-3-na-ldr-primal} with $ t-1 $ (when $t=1$ there is no second term), we can replace $ \Exp_{\xiT}\left[\Phi_{tk}(\xisupt)^\top B_t(\xisupt)\Exp_{\xiprimeT}[y^{\star}_{t-1}(\xiprimeT) | \xiprimesuptminusone = \xisuptminusone] \right] $ with $
	\Exp\left[\Phi_{tk}(\xisupt)^\top B_t(\xisupt)y^{\star}_{t-1}(\xiT) \right] $  \rev{since by construction there exists} $ \kprimeindex\in\Kprimet $ with 
	$ \Psi_{{t-1}\kprimeindex}(\xiT)^\top =\Phi_{tk}(\xisupt)^\top B_t(\xisupt) $. Putting all three terms of \eqref{eq:break}
	back together, we obtain that the expression \eqref{eq:lhs} is equal to
	\begin{equation*} 
	\Exp_{\xiT}\left[\Phi_{tk}(\xisupt)^\top \left(A_t(\xisupt)y^{\star}_t(\xiT) + B_t(\xisupt)y^{\star}_{t-1}(\xiT) - b_t(\xisupt)\right)\right].
	\end{equation*}
	 Since $ y_t(\xiT)\in Y_t(\xiT) $, this expectation is equal to zero and the constraint \eqref{eq:const-strength-ldr-state} is satisfied. \Halmos \endproof

Thus, we conclude that if the basis functions are selected carefully, then the restricted NA dual is not worse than the restricted SW dual. Our numerical experiments in Section \ref{sec:results} illustrate that it can indeed provide strictly better bounds.


\subsection{\rev{Choice of Basis Functions}}\label{sec:basisfunctions}
\rev{
Basis function selection appears as a design choice in many computational fields, and as such has a rich literature. As mentioned by \cite{powell2009you} in the context of approximate dynamic programming, this choice is problem-specific and requires the knowledge of underlying features that affect the future for a proper design. A natural idea is to define one basis function per uncertain parameter, referred to as the \emph{standard basis functions}. The following proposition shows that, when using the standard basis functions in the NA dual, the basis functions used in stage $t$ can be restricted to use only the uncertain parameters for periods from stage $t$ to the end of the horizon.

\begin{proposition}\label{prop:nohistory}
    Let $\Psi_t(\xiT) = (1, \xiparrand_2, \xiparrand_3, \dots, \xiparrand_T)$. There exists an optimal solution $\naldrvec^\star$ to problem \eqref{eq: ldr-na-lagrange-dual} where $\naldrvec_{t\kprimeindex}^\star = 0$ for $1 < \kprimeindex \leq t $.
\end{proposition}
\proof{Proof}
For $\kprimeindex \leq t$, $\Exp_{\xiprimeT}[\ \xiparrand'_{\kprimeindex}\  | \xiprimesupt = \xisupt] = \xiparrand_{\kprimeindex}$. Hence the coefficient of $\naldrvec_{t\kprimeindex}$ for $1 < \kprimeindex \leq t $ in the objective function of \eqref{eq:na_subpproblem} is 0 and we can fix $\naldrvec_{t \kprimeindex}^\star = 0$ for $1 < \kprimeindex \leq t $ in an optimal solution.
\Halmos \endproof

The following proposition confirms the intuition that including a basis function that is a linear combination of other basis functions would not change the bound obtained.

\begin{proposition}\label{prop:nolinearcomb}
    For $t\in[T]$, let $\Psi^1_t=(\Psi_{t1},\dots,\Psi_{t \Kprimet})$ be a given set of basis functions  and for $w^t \in \mathbb{R}^{\Kprimet}$ define $\Psi^2_t=(\Psi_{t1},\dots,\Psi_{t \Kprimet},(w^t)^\top \Psi^1_t)$. Denoting by $\nu_R^{\text{NA}}(\Psi)$ the NA dual bound obtained by using basis functions $\Psi$, we have $$\nu_R^{\text{NA}}(\Psi^1) = \nu_R^{\text{NA}}(\Psi^2).$$
\end{proposition}
\proof{Proof}
In the NA primal characterization, the choice of basis functions affects constraints \eqref{eq:const-3-na-ldr-primal}. For the basis function $(w^t)^\top\Psi^1_t$ we have $
\Exp_{\xiT}\Big[ \big((w^t)^\top\Psi^1_{t}(\xiT)\big)^\top \big(y_t(\xiT) - \Exp_{\xiprimeT}[y_t(\xiprimeT) | \xiprimesupt = \xisupt]\big) \Big] = \sum_{k\in \Kprimet}w^t_{k} \Exp_{\xiT}\Big[\Psi_{tk}(\xiT)^\top \big(y_t(\xiT) - \Exp_{\xiprimeT}[y_t(\xiprimeT) | \xiprimesupt = \xisupt]\big) \Big] = 0$, which is a linear combination of the constraints associated with the basis functions included in $\Psi^1_t$, and thus is redundant. 
\Halmos \endproof

If we have some knowledge that a certain function of the uncertain parameters impacts the model, it may be useful to use this function as a basis function. For instance, in the inventory example studied by \cite{bodur2018two}, demand at stage $t$ is obtained from the multiplication of random variables observed in the prior stages. Therefore, a basis function at stage $t$ is defined by multiplying those random variables (e.g., as opposed to using the individual random variables as basis functions). For a more complex stochastic process, the literature of approximate dynamic programming, approximate linear programming and machine learning offer more nonlinear choices, e.g., polynomial, piecewise constant/linear \citep{friedman2001elements} and radial basis functions \citep{buhmann2000radial}. Rather than fixing the form of the basis functions, some works have opted for  design procedures that leverage learning methods such as neighbourhood component analysis      \citep{keller2006automatic} and neural networks  \citep{hajizadeh2018optimized} for constructing basis functions. In Sections \ref{sec:bf-selection} and \ref{sec:bounds-RWA} we study the impact of different options of basis functions for our two problem classes and show the importance of their proper design on the solution time and quality.
}


\subsection{Solving the Restricted Dual Problems}\label{sec:SAA}

The restricted dual problems \eqref{eq: ldr-lagrange-dual} and \eqref{eq: ldr-na-lagrange-dual} have the form of maximizing an expected value of a nonsmooth concave function of the decision variables (the coefficients of the LDDR policies). A wide variety of algorithms for (approximately) solving such problems exist, including stochastic approximation based methods (e.g., \cite{robbinsmonro51,robustsa09}), stochastic decomposition \citep{stochdecomp}, and sample average approximation (SAA) (e.g., \cite{shapiro2009lectures}). For concreteness, we describe an SAA approach.  Let $\{ \boldsymbol{\xi}^T_\omega \}_{\omega \in \Omega}$ be a given sample of sample-path scenarios and $p_\omega$ denote the probability of scenario $\omega \in \Omega$. For example, if the scenarios are generated via Monte Carlo sampling, then $ p_{\omega} = 1/|\Omega| $ for all $\omega \in \Omega$. SAA replaces the expectations in \eqref{eq: ldr-lagrange-dual} and \eqref{eq: ldr-na-lagrange-dual} with sample averages, which respectively give the following SAA models for approximating the restricted SW and NA duals
\bsub
\label{eq:SAA-sw}
\begin{alignat}{2} 
\max_{ \theta, \{\ldrvec_t\}_{t \in \set{T}}}  \   &
\sum_{\omega \in \Omega} p_\omega \theta_{\omega} 
\\
\text{s.t.} \ \ & \theta_{\omega} \leq \sum_{t \in [T]}
\mathcal{L}^{\text{SW}}_t(\ldrvec,\boldsymbol{\xi}_\omega^t), \qquad && \omega \in \Omega \label{eq:SAA-sw-bcuts}
\end{alignat}
\esub
and
\bsub
\label{eq:SAA-NA}
\begin{alignat}{2}
\max_{ \eta, \{ \naldrvec_t \}_{t \in \set{T}}}  \   &
\sum_{\omega \in \Omega} p_\omega \eta_{\omega} 
\\
\text{s.t.} \ \ & \eta_{\omega} \leq \mathcal{L}^{\text{NA}}(\naldrvec,\boldsymbol{\xi}^T_\omega), \qquad && \omega \in
\Omega. \label{eq:SAA-na-bcuts} 
\end{alignat}
\esub
In order to solve \eqref{eq:SAA-sw} and \eqref{eq:SAA-NA}, we use a regularized Benders method
\citep{ruszczynski1986regularized}. It is proven that under certain conditions (e.g., the feasible sets in \eqref{eq: ldr-lagrange-dual} and \eqref{eq: ldr-na-lagrange-dual} being nonempty and bounded, and the  expectations being finite), as the sample size increases, the solutions to the SAA models converge to the ones of the respective original  problems \citep{shapiro2009lectures}. To assure convergence of the SAA problems, it may be required to put bounds on the LDDR decision variables.

We refer the reader to \citep{kiwiel95,ruszczynski1986regularized} for details and implementation strategies for the
regularized Benders method. This algorithm uses a master problem that approximates the constraints \eqref{eq:SAA-sw-bcuts}
and \eqref{eq:SAA-na-bcuts}, respectively, with a finite set of Benders cuts. 
The part of the algorithm that requires specialization for its application to problems \eqref{eq:SAA-sw} and \eqref{eq:SAA-NA} is the specification of the Benders cuts that are added to the master problem given a current master problem solution $\hat{\ldrval}$ or $\hat{\naldrval}$, for problems  \eqref{eq:SAA-sw} and \eqref{eq:SAA-NA}, respectively.

 For problem \eqref{eq:SAA-sw},  subproblems evaluating $\mathcal{L}^{\text{SW}}_t(\hat{\ldrvec},\boldsymbol{\xi}^t_{\omega})$ are solved for $\omega \in \Omega$ and $t\in\set{T}$. Letting $\hat{x}_{t}^{\omega}$ denote the optimal solution of subproblem for $t$ and scenario $\omega$, the Benders cut  is given below
\begin{align*}
\theta_{\omega} \leq \sum_{t \in \set{T}} \mathcal{L}^{\text{SW}}_t(\hat{\ldrvec},\boldsymbol{\xi}_\omega^t) & +
\sum_{t \in \set{T}} \sum_{k \in [\Kt]} (\ldrval_{tk} - \hat{\ldrval}_{tk}) \, g_{tk}(\hat{\ldrvec},\boldsymbol{\xi}_\omega^t)  
\end{align*}
where $ g_{tk}(\hat{\ldrvec},\boldsymbol{\xi}_\omega^t) = -\Phi_{tk}(\boldsymbol{\xi}_\omega^t)b_t(\boldsymbol{\xi}_\omega^t) + \left(\Phi_{tk}(\boldsymbol{\xi}_\omega^t)A_t(\boldsymbol{\xi}_\omega^t) + B_{t}(\boldsymbol{\xi}_\omega^t)\Exp_{\xiprimeT}\left[\ \Phi_{tk}(\xiprimesupt)\ \ |  \ \boldsymbol{\xi}_\omega^{t-1} \right]\right)^\top \hat{x}^\omega_t $ 
is the $k^\text{th}$ component of the subgradient of $\mathcal{L}^{\text{SW}}_t(\ldrvec,\boldsymbol{\xi}_\omega^t)$ with respect to $ \ldrvec $, at point $ \hat{\ldrvec}$ evaluated at $\boldsymbol{\xi}_\omega$.  

For problem \eqref{eq:SAA-NA}, subproblems evaluating $\mathcal{L}^{\text{NA}}(\hat{\naldrvec},\boldsymbol{\xi}_{\omega}^T)$ for each $\omega \in \Omega$ are solved.
Denoting the subproblem optimal solutions by $\hat{y}^\omega$ for each $\omega \in \Omega$, the Benders optimality cut for the restricted NA Lagrangian dual is the inequality
\begin{align*}
\eta_{\omega} \leq & \  \mathcal{L}^{\text{NA}}(\hat{\naldrvec},\boldsymbol{\xi}^T_\omega)
+ \sum_{t \in \set{T}} \sum \limits_{\kprimeindex \in [\Kprimet]} \bigg( \Psi_{t \kprimeindex}(\boldsymbol{\xi}^T_\omega) - \Exp_{\xiprimeT}[\Psi_{t \kprimeindex}(\xiprimeT) | \boldsymbol{\xi}_\omega^t] \bigg) 
\Big[ 
\>  (\naldrval_{t \kprimeindex}-\hat{\naldrval}_{t \kprimeindex})   \hat{y}_{t}^{\omega}
\Big].
\end{align*}
Note that in the presence of stagewise independence or a recursive form such as an autoregressive process, $ \Exp_{\xiprimeT}[\Psi_{t \kprimeindex}(\xiprimeT) | \boldsymbol{\xi}_\omega^t]  $ can be computed directly. Otherwise conditional expectations need to be approximated, e.g., by sampling.

Solving \eqref{eq:SAA-sw} and \eqref{eq:SAA-NA} yields candidate LDDR solutions, say $ \ldrvec^\star_{\Omega} $ and $ \naldrvec^\star_{\Omega} $, respectively. Let $\{ \boldsymbol{\xi}^T_\omega \}_{\omega \in \Omega'}$ be an independent evaluation sample with $|\Omega'| >> |\Omega|$. In order to get statistically valid lower bounds on the optimal value of the original {\MSMIP} problem, subproblems $ \mathcal{L}^{\text{SW}}_t(\cdot) $ and $ \mathcal{L}^{\text{NA}}(\cdot) $ respectively given in \eqref{eq:SAA-sw} and \eqref{eq:SAA-NA},  are solved with fixed $ \ldrvec^\star_{\Omega} $ and $ \naldrvec^\star_{\Omega} $ for every scenario in the evaluation sample, i.e., solving $ \mathcal{L}^{\text{SW}}_t(\ldrvec^\star_{\Omega}, \boldsymbol{\xi}_{\omega}^t) $ and $ \mathcal{L}^{\text{NA}}(\naldrvec^\star_{\Omega}, \boldsymbol{\xi}_{\omega}) $ for all $ \omega\in\Omega' $. The lower end of a confidence interval based on the obtained values is a statistically valid lower bound, regardless of how $\ldrvec^\star_{\Omega} $ and $ \naldrvec^\star_{\Omega} $ were obtained (e.g., these need not be optimal solutions).

A summary of the notation used in the two proposed restricted Lagrangian duals is given in Table \ref{tab:duals-compare}. 
\begin{table}[tbp]
    \renewcommand{\arraystretch}{0.55}
	\centering
	\begin{tabular}{llrl}
		\toprule
		& Stagewise dual (SW) &       & Nonanticipative dual (NA)\\
		\midrule
		Relaxation & State equations &       & Nonanticipativity constraints \\
		Dual variables & $ \pi_t $ &       & $ \gamma_t $  \\
		LDDR decision variables & $ \ldrvec_t$ &       & $ \naldrvec_t$ \\
		Basis Functions & $ \Phi_t $ & & $ \Psi_t $\\
		Subproblem per & $ t $ and $ \xisupt $ &       & $\xiT$ \\
		Optimal value & $ \nu^{\text{SW}} $ & \multicolumn{1}{l}{$\leq$} &  $\nu^{\text{NA}} $\\
		\bottomrule
	\end{tabular}%
	\caption{Summary of the characteristics of the SW and NA duals}
	\label{tab:duals-compare}
\end{table}%
\spacing


\section{Primal Policies}\label{sec:ub-policy}
We next describe how restricted dual solutions can be used to obtain primal policies. \rev{We present a classical policy based on solving deterministic approximations, and two policies that use solutions from the restricted duals.

All the primal policies follow a \revtwo{folding-horizon} horizon framework: at every stage we make a decision respecting the current stage constraints, then after observing the value of the random parameters, we \emph{\revtwo{fold} the horizon}, meaning that the decisions of the current stage are passed as parameters to the next stage and new decisions are made. More formally, let $\{ \xisupomegaval \}_{\omega \in \Omega'}$ be a collection of sample path observations (scenarios) that we will use to evaluate a policy. For each scenario $ \xisupomegaval $,  the cost of using this policy for this scenario is computed as $ U_\omega = \sum_{t\in\set{T}} c_t(\xisuptomegaval)^\top \hat{x}_t(\xisuptomegaval)$, where $ \hat{x}_t(\xisuptomegaval), t \in \set{T}$ are the decisions chosen by the policy. In each of the policies we consider, the solution $\hat{x}_t(\xisuptomegaval)$ is obtained by solving an optimization model that takes as input $\hat{x}_{t-1}(\xisuptmoneomegaval)$ (when $t>1$) and the current random parameter observations $\xisuptomegaval$. The optimization model solved at stage $t$ includes, at a minimum, the constraints:
\begin{align*}
& A_t(\xisuptomegaval)x_t  = \rev{b_t(\xisuptomegaval)} - B_t(\xisuptomegaval)\hat{x}_{t-1}(\xisuptmoneomegaval), \\
& x_{t} \in X_t(\xisuptomegaval). 
\end{align*}
After solving the optimization problem at stage $t$, $\hat{x}_t(\xisuptomegaval)$ is set to the $x_t$ component of the optimal solution. These constraints, together with the relatively complete recourse assumption (Assumption \ref{ass:rel-recourse}) ensure that a feasible solution exists at every stage $t$, and hence each of the presented policies is feasible. A confidence interval of the expected cost of each policy can then be computed using the values in $\{ U_\omega\}_{\omega\in\Omega'} $, since each of these values represents an independent observation of the cost of the simulated policy. The upper end of this confidence interval provides a statistical upper bound on the cost of the associated feasible policy, and therefore an upper bound on the optimal value of the \MSMIP \ model. 

In the following sections we discuss three primal policies that share the same \revtwo{folding-horizon} framework and differ from each other in how the decisions of the current stage are made, i.e., the optimization problems that need to be solved. More specifically, the key difference in these approaches is in how they approximate the impact of the current decisions on future costs.} 


\subsection{Classical Approach}
We first review a classical approach which obtains a policy by replacing the uncertain future parameter values with their conditional expected value. For a partial sample path $\xisupt$, denote by $ \xiaverageval$ the conditional expected value of $\xi^s$ given $\xisupt$, for $s \geq t$. \rev{Within the \revtwo{folding-horizon} framework, at stage $t$, given the solution of the previous stage $\hat{x}_{t-1}(\xisuptmoneomegaval)$, we solve the following \emph{deterministic} problem (initialized at $t=1$ with $\hat{x}_1$ as the vector of zeros):}
\bsub
\label{eqs:2SPext-BasePrimal}
\begin{alignat}{2}
\min_{ \revtwo{\{x_t \}_{t \in \set{t,T}} }} \ \ & c_t(\xisuptomegaval)^\top x_t +  \sum_{s \in [t+1,T]} c_s(\xiaverageval)^\top x_{s}  \label{eq:obj-2SPext-BasePrimal} \\
\text{s.t.} \ \ & A_t(\xisuptomegaval)x_t  = \rev{b_t(\xisuptomegaval)} - B_t(\xisuptomegaval)\hat{x}_{t-1}(\xisuptmoneomegaval),  \label{eq:const1-2SPext-BasePrimal} \\
& A_s(\xiaverageval)\ x_s  + B_s(\xiaverageval)\ x_{s-1}= b_s(\xiaverageval), && s \in [t+1,T] \label{eq:const4-2SPext-BasePrimal} \\
& x_{s} \in X_s(\xiaverageval), &\quad& s \in [t,T].  \label{eq:const5-2SPext-BasePrimal} 
\end{alignat}
\esub
We refer to the policy obtained from this upper bounding procedure as the \emph{conditional expected value policy}.
\rev{
Note that, for a sample path $\xisupomegaval$, only the state equations of the current stage are enforced. For future stages, the state equations are only written for $\xiaverageval$ and not almost every scenario. Therefore, $\hat{x}_t(\xisuptomegaval)$ might not satisfy the state equations in the next stage for $\xiparval^{t+1}_{\omega}$. In the next section, we take advantage of our LDDR solutions with the hope to identify solutions that provide lower cost over the full horizon.
}


\subsection{Restricted Stagewise Lagrangian Dual Driven Policy}
To improve the conditional expected value policy, an idea is to add a penalty obtained from the SW Lagrangian dual 
\rev{
for violating the state equations at the next immediate stage for $\xiparval^{t+1}_{\omega}$.  
}
Given a restricted SW dual solution $ \hat{\pi}_t(\xisupt) =  \Phi_{t}(\xisupt) \ \hat{\beta}_t $, $ t\in\set{T}$, we obtain the \emph{SW dual driven policy} by following the same \revtwo{folding-horizon} procedure above 
\rev{
with the objective function as the convex combination of \eqref{eq:obj-2SPext-BasePrimal} and the coefficient of $x_t$ decision variables in $\mathcal{L}^{\text{SW}}_t$: 
\bsub
\begin{alignat}{2}
\min_{ \revtwo{\{x_s \}_{s \in \set{t,T}} }} \ \ & \Big(c_t(\xisuptomegaval) + \lambda\ \big(\hat{\ldrvec}^\top_{t}\Phi^\top_{t}(\xisuptomegaval)\ A_t(\xisuptomegaval) + \hat{\beta}^\top_{t+1} \ \Exp_{\xiprimeT}\big[ \Phi^\top_{t+1}(\xiprimesuptplusone) \  B_{t+1}(\xiprimesuptplusone) \ | \ \xisuptomegaval \big]\big)\Big)^\top  x_t \label{eq:modified-exp-up-obj} \\
& + (1-\lambda) \sum_{s \in [t+1,T]} c_s(\xiaverageval)^\top x_{s}  \nonumber \\
\text{s.t.} \ \ & A_t(\xisuptomegaval)x_t  = \rev{b_t(\xisuptomegaval)} - B_t(\xisuptomegaval)\hat{x}_{t-1}(\xisuptmoneomegaval),  \nonumber \\
& A_s(\xiaverageval)\ x_s  + B_s(\xiaverageval)\ x_{s-1}= b_s(\xiaverageval), && s \in [t+1,T] \nonumber \\
& x_{s} \in X_s(\xiaverageval), &\quad& s \in [t,T].  \nonumber
\end{alignat}
\esub
where $\lambda\in [0,1]  $ is 
the weight of information from the SW Lagrangian dual and 
a parameter of the policy. Setting $\lambda = 0$ leads to the conditional expected value policy, while $\lambda = 1$ corresponds to a modified myopic policy where the only information about the future that we pass to the primal algorithm is the penalty associated with the next immediate stage.}


\subsection{Restricted Nonanticipative Lagrangian Dual Driven Policy}
Given a restricted NA dual policy $ \hat{\gamma}_t(\xiT) = \ \Psi_t(\xiT)\hat{\alpha}_t$, for all $ t\in\set{T} $, we define the \emph{NA dual-driven policy} as follows. \rev{Using the discussed \revtwo{folding-horizon} framework,} at stage $ t $, with observed history $ \xisuptomegaval $ and previous stage decisions $ \hat{x}_{t-1}(\xisuptmoneomegaval) $, choose $ x_t(\xisuptomegaval)$ as an approximate solution of the following \emph{two-stage stochastic} program
\bsub
 \label{eq:twostage-master}
\begin{alignat}{2}    
\min_{x_t}  \ \ & c_t(\xisuptomegaval)^\top x_t + \Exp_{\xiT}\left[g_t(x_t,\xiT)\ |  \ \xisuptomegaval \right] \label{eq:twostage-master-obj} \\
\text{s.t.} \ \ & A_t(\xisuptomegaval)x_t  = b_t(\xisuptomegaval) - B_t(\xisuptomegaval)\hat{x}_{t-1} \\
& x_t\in X_t(\xisuptomegaval) 
\end{alignat}
\esub
where
\bsub
\label{eq:na-primal-inner}
\begin{alignat}{2}
g_t(x_t,\xiT) = \min_{ \{x_s\}_{s \in [t+1,T]}} \ \ & \sum_{s \in [t+1,T]} \left( c_s(\xisups) + \hat{\gamma}_s(\xiT)  - \Exp_{\xiprimeT}[\hat{\gamma}_s(\xiprimeT) \ | \ \xisupt]  \right)^\top x_s \label{eq:obj-primal-na_sub} \\
\text{s.t.} \ \ & A_s(\xisups)x_s + B_s(\xisupsminusone)x_{s-1} = b_s(\xisups), &\quad& s \in [t+1,T]  \label{eq:const1-primal-na_sub} \\
& x_s\in X_s(\xisups), && s \in [t+1,T].  \label{eq:const2-primal-na_sub} 
\end{alignat}
\esub

In order to solve the problem \eqref{eq:twostage-master}, for each realization $ \xisuptomegaval $, a sample $ \{\xitwostagevalomega\}_{\omega\in\Omega_{|t}} $ is generated, where
$\xitwostagerand$ is a random variable representing the (conditional) scenarios after the history up to stage $ t $ has been observed. Then we solve it by replacing the expectation in \eqref{eq:twostage-master-obj} with a sample average using this sample. 
\rev{Since this scheme is used in a \revtwo{folding-horizon} fashion,} in order to estimate the expected cost of this policy, multiple sample paths $ \xi^T $ are generated, and for each one a cost is recorded by the scheme outlined above. Finally, a confidence interval on upper bound value is obtained over these costs.

As the second stage problem \eqref{eq:na-primal-inner} is a MIP, given a set of conditional scenarios for the future, we can solve the extensive form of this two-stage problem, by creating a copy of the second-stage variables in problem \eqref{eq:na-primal-inner} for each conditional scenario and embedding them in \eqref{eq:twostage-master} for all scenarios.
\rev{
Note that although two-stage stochastic MIPs are computationally difficult to solve, there has been significant research into methods for solving this problem class. This includes the integer L-shaped method \citep{laporte1993integer} and a large literature on methods for generating improved cuts within such a solution approach, e.g., \citep{angulo2016improving,bodurluedtke:17,qi2017ancestral,ntaimotannerdisj:08,gadesimgesen:12,senhigle:05,sheralizhu:07,sensherali:mp06,minjiaosimge:14,bendersdual:20,van2021converging}. Another general approach is dual decomposition \citep{caroe1999dual}, see \citep{kim2019} for extensions and a recent parallel implementation.
\revtwo{Approximation methods for solving two-stage stochastic MIPs have also been proposed, e.g., \citep{sandikcci2013hierarchy,romeijnders2016convex,van2021loose}}.
See \citep{kuccukyavuz2017introduction} for a more extensive overview of exact and approximate solution methods for two-stage stochastic MIPs. 
}

Although two-stage stochastic MIP problems are  challenging, when implementing this policy it needs to only be solved once per time stage.
On the other hand, estimating the expected value of this policy is computationally challenging since it is required to simulate many sample paths to construct a confidence interval. Thus, it may be necessary to heuristically limit the effort in solving the two-stage stochastic MIP, e.g., by using a small number of scenarios, or by terminating the solution process after the solver processes just a limited number of branch-and-bound nodes. Another approach is to generate a larger number of scenarios and then representing them by a smaller size sample, e.g., using a scenario reduction technique like clustering.

\begin{figure}[htbp]
	\centering
	\includegraphics[scale=0.88]{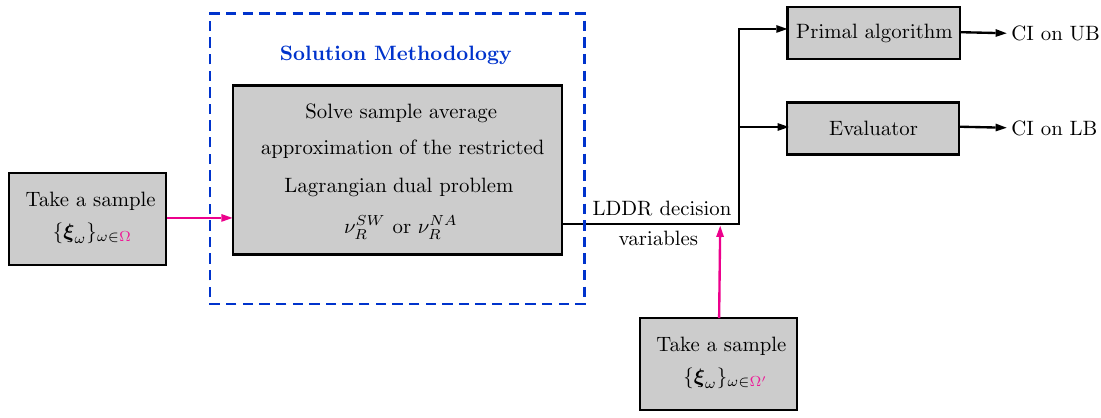}
	\caption{Solution framework}
	\label{fig:framework}
\end{figure}

The overall framework for providing bounds on {\MSMIP} is illustrated in Figure \ref{fig:framework}. Using LDDRs, one
can obtain a Lagrangian dual policy. The parameters of this policy are then passed to an evaluator to provide confidence intervals (CIs) on the out of sample scenarios. When comparing upper bounds obtained from different methods, the same set of scenarios is used in their evaluation.


\section{Computational Experiments}\label{sec:results}
We illustrate our proposed approach on \rev{two problems: ($i$) multi-item lot-sizing problem with backlogging and production lag (MSLot), and ($ii$) routing and wavelength assignment problem (MSRWA)}. We present our analysis on the choice of basis functions, as well as comparison between the bounds of various proposed methods.


\subsection{Implementation Details}\label{sec:implementation}
To solve the SAA problems \eqref{eq:SAA-sw} and \eqref{eq:SAA-NA} for the MSLot and MSRWA, we have used the regularized Benders method, with all the parameters of the algorithm being exactly as provided by \cite{lubin2013parallelizing}. The convergence tolerance of the algorithm is set to 0.001. 
The algorithms are implemented in C++ using IBM ILOG CPLEX 12.8 for solving MIPs, and the experiments are conducted on a MacOS X with 2.3 GHz Intel Core i5 and 16 GB RAM. \rev{For all bounding methods we have considered a maximum running time of 6 hours.}

For \revtwo{ dual bounds, we use the restricted SW and NA dual problems. For primal bounds, we use the dual driven primal policies. The SW dual driven policy is implemented with $ \lambda = 0.25 $ in the objective function \eqref{eq:modified-exp-up-obj}}. \rev{This value is selected based on some preliminary analysis on the \revtwo{lot-sizing} problem (whose results are provided in Table \ref{tab:diff-lambda} of the Appendix).} For \revtwo{the NA dual driven primal policy,} the two-stage stochastic MIP problem \eqref{eq:twostage-master} is approximately
solved using a sample with $ |\Omega_{|t}| =25 $. This set is generated by randomly sampling 100 scenarios, then
clustering these into 24 groups and using the mean of each group as a scenario in the sample, and adding the conditional
expectation $ \xiaverageval $ as the last scenario in the sample.


\subsection{Benchmarks}\label{sec:benchmarks}
The same set of scenarios is used for estimating lower bounds and evaluating primal policies. The perfect information (PI) bound is used as our benchmark for lower bounds, which is obtained similar to the restricted NA lower bound, with the LDDR
variables, $ \naldrvec $, fixed to 0. This is the bound  obtained by solving the version of the problem where the complete information about the future is assumed to be known. We use the conditional expected value policy as the benchmark upper bounding policy. \rev{Since this formulation is not based on a scenario-tree model (and approximating it with a scenario tree would require an intractable number of nodes), to the best of our knowledge, these are the best available comparison methods for the general \MSMIP \ problem class.} 


\subsection{Multi-item Stochastic Lot-sizing Problem}
The {\MSMIP} formulation for the MSLot problem is as follows
\bsub
\label{eq:lot-sizing}
\begin{alignat}{2}
\hspace{-0.1cm}\min  \ \, &  \Exp\Bigg[\sum\limits_{t \in \set{T}} \bigg( \sum\limits_{j \in \set{J}} \Big(
C_{tj}^{i^+}(\xisupt)i^+_{tj}(\xisupt) + C_{tj}^{i^-}(\xisupt)i^-_{tj}(\xisupt) + C_{tj}^{y}(\xisupt)y_{tj}(\xisupt) \Big) +   C_{t}^{o}(\xisupt)o_{t}(\xisupt)\bigg)
\Bigg] \label{eq:lot-sizing-obj}  \\
\hspace{-0.1cm}\text{s.t.} \ & i^-_{tj}(\xisupt)-i^+_{tj}(\xisupt) + i^+_{t-1,j}(\xisuptminusone) - i^-_{t-1,j}(\xisuptminusone) + x_{t-1,j}(\xisuptminusone) = D_{tj}(\xisupt),   &&  \nonumber
\\
&&&\hspace{-36mm}t \in \set{T}, j \in \set{J}, \allxi \label{eq:lot-sizing-state-equations}\\
\hspace{-0.1cm}&  \sum\limits_{j \in \set{J}}(TS_jy_{tj}(\xisupt) +TB_jx_{tj}(\xisupt)) - o_{t}(\xisupt)\leq C_t, & \, & \hspace{-36mm}t \in \set{T}, \allxi \label{eq:lot-sizing-overtime}\\
\hspace{-0.1cm}& M_{tj} \, y_{tj}(\xisupt) - x_{tj}(\xisupt) \geq 0, && \hspace{-36mm}t \in \set{T}, j \in \set{J}, \allxi  \label{eq:lot-sizing-setup}\\
\hspace{-0.1cm}& i^+_{tj}(\xisupt) \leq I_{tj}, && \hspace{-36mm}t \in \set{T}, j \in \set{J}, \allxi \\
\hspace{-0.1cm}&i^+_{tj}(\xisupt) + x_{tj}(\xisupt) \leq I_{t+1,j}, && \hspace{-36mm}t \in \set{T}, j \in \set{J}, \allxi  \\
\hspace{-0.1cm}& 0 \leq o_{t}(\xisupt) \leq O_{t}, && \hspace{-36mm}t \in \set{T}, \allxi \\
\hspace{-0.1cm}& x_{tj}(\xisupt),i^+_{tj}(\xisupt), i^-_{tj}(\xisupt) \geq 0, && \hspace{-36mm}t \in \set{T}, j \in \set{J}, \allxi  \\
\hspace{-0.1cm}& y_{tj}(\xisupt) \in\{0,1\}, && \hspace{-36mm}t \in \set{T}, j \in \set{J}, \allxi \label{eq:lot-sizing-integrality}
\end{alignat}
\esub
where $ x_{tj}(\xisupt), i^+_{tj}(\xisupt), i^-_{tj}(\xisupt) $ are decision variables representing production level, inventory and backlog of product type $ j $ at stage $ t $, respectively. Binary decision variable $ y_{tj}(\xisupt) $ is  equal to 1 if production of item $ j $ is setup at stage $ t $, 0 otherwise. Decision variable	$ o_{t}(\xisupt) $ measures the overtime at stage $ t $. $ C_{t}, I_{tj}, O_{t}$ are the production capacity, inventory capacity of product $ j $, and the overtime bound at stage $ t $, respectively. $ D_{tj}(\xisupt) $ is the demand of product $ j $ at stage $ t $.
$ C^{i^+}_{tj}(\xisupt), C^{i^-}_{tj}(\xisupt), C^{y}_{tj}(\xisupt)  $ are respectively the costs of  holding, backlog and fixed setup for product $ j $, and $ C^{o}_{t}(\xisupt) $ is the overtime cost at stage $ t $. For product $ j $, $ TS_j $ is the setup time, while
$ TB_j $ is the production time per unit. Constraints \eqref{eq:lot-sizing-state-equations} are the state equations,
linking the inventory, backlog and production of consecutive stages. Note that there is a production lag of 1, meaning
that the amount that is produced at stage $ t $ is not available until stage $ t+1 $. Overtime is measured by constraints \eqref{eq:lot-sizing-overtime}, while the same set of constraints ensure the production capacity is respected. Constraints \eqref{eq:lot-sizing-setup} link the production and setup decisions,  using sufficiently large big-$M$ values in the absence of a given limit on the production quantity. The rest of the constraints determine the bounds and integrality constraints on the decision variables. The objective function \eqref{eq:lot-sizing-obj} is the total expected cost, including the costs of holding, backlog, setup and overtime.


\subsubsection{Base Data}\label{sec:data-instances}
The dataset of our experiments is loosely based on the work of  \cite{helber2013dynamic}. 
We consider the following autoregressive process model for representing the correlation between demands of different stages
\begin{equation*}
Y_{t+1,j} = \rho \, Y_{tj} + (1-\rho) \, \epsilon_{t+1,j},
\end{equation*}
where $ \epsilon_{tj} $ is a lognormal random variable with mean 1 and standard deviation 0.5. Demands are modeled as 
$$ D_{tj} = \rho^Y Y_{tj} \, \mu_{tj} + (1-\rho^Y) \, \delta_{tj},$$ 
where $ \mu_{tj}$ is the mean of the demands of product type $ j $ at stage $ t $, and $ \delta_{tj} $ is a  random variable having a mean of $
\mu_{tj} $. Thus, in this model there is an underlying autoregressive process ($ Y_{tj} $) and the demand in each period is partially driven by this and also by an external random variable. This assures that demands are non-negative. We consider a lognormal distribution for $ \delta_{tj} $ whose standard deviation is  $0.2 \, t \, \mu_{tj}$, reflecting higher demand uncertainty further in the future. For $  \mu_{tj} $ values, we use the means of the demands provided by \cite{helber2013dynamic}.
\revtwo{The rest of the configurations can be found in 
the Appendix.}


\subsubsection{Basis Function Selection}\label{sec:bf-selection}
As mentioned in \rev{Section \ref{sec:basisfunctions}}, the choice of basis functions is not restricted by any predetermined form. Therefore, in choosing the basis functions for the MSLot problem we have options. 
\revtwo{At stage $ \hat{t} $, for the SW dual, using the standard basis functions means using $ D_{tj} $ for all $ j $ and $ t \leq \hat{t} $, while in the NA dual it involves the complete set of demands in all stages of the planning horizon for all products. Table \ref{tab:options} summarizes some of the alternatives  considered in our experiments.}
Note that, Option 1 in the NA dual does not span over the entire planning horizon and it ignores the observed history. \rev{According to Proposition \ref{prop:nohistory}, these} two are equivalent. 

\begin{table}[htbp]
\renewcommand{\arraystretch}{0.55}
 \centering
 \caption{Various options to be used in basis functions, at stage $ \hat{t} $ for product $ \hat{j} $}
    \begin{tabular}{c cll l cll l}
    \toprule
    Option & & \multicolumn{3}{c}{SW dual} & & \multicolumn{3}{c}{NA dual} \\
    \midrule
    1     & & $1 \ ; D_{tj}$, & $\forall j$, & $\forall t\leq \hat{t}$ &&  $1\ ; D_{tj}$, & $\forall j$, & $\forall t\leq T: \hat{t} < t$ \\
    2     & & $1\ ; D_{tj}$, & $\forall j$, & $t=\hat{t}$ & & $1\ ; D_{tj}$, & $\forall j$, & $t=\hat{t}+1$ \\
    3     & & $1\ ; D_{tj}$, & $j=\hat{j}$, & $\forall t\leq \hat{t}$ & & $1\ ; D_{tj}$, & $j=\hat{j}$,& $ \forall t\leq T: \hat{t} < t $ \\
    4     & & $1\ ; D_{tj}$, & $j=\hat{j}$,  & $t= \hat{t}$ & & $1\ ; D_{tj}$, & $j=\hat{j}$, & $ t = \hat{t}+1$ \\
    \bottomrule
    \end{tabular}%
 \label{tab:options}%
\end{table}%

In Table \ref{tab:options}, Options 2 and 4 in  the SW dual only use the stage $\hat{t}$ demands as basis functions for
stage $\hat{t}$ constraints, and in the case of the NA dual, they only use the stage $\hat{t}+1$ demands for basis
functions used to relax stage $\hat{t}$ constraints. Options 3 and 4 only use product $\hat{j}$ demands as basis functions for constraints associated with product $\hat{j}$. Thus, Option 1 has the most basis functions, and hence should have the best bound, whereas Option 4 uses the fewest basis functions. We compare these four options to determine which one gives the best trade-off between the quality of the bound and the computational effort. 
\rev{For instances with $J=3, \rho = 0.6, \rho^Y = 0.2$,} Tables \ref{tab:bf-options-sw} and \ref{tab:bf-options-na} present the solution time in seconds and the lower bounds returned by the SW and NA duals, respectively. For this comparison, three instances with $ T=4, 6, 8 $ are solved using the four basis function options \rev{and only penalizing the NA constraints on $x$ variables (see Section \ref{sec:na-var-selection})}. The bounds are the means of the confidence intervals over the objective values for all the scenarios in the evaluation sample. They are scaled such that 100 is the best known bound obtained for that instance. \rev{The best bounds for these instances are obtained by an option of the NA dual reported in Table \ref{tab:na-var-selection}, which is why all the scaled bounds in Tables \ref{tab:bf-options-sw} and \ref{tab:bf-options-na} are less than 100.}
\begin{table}[htbp]
\renewcommand{\arraystretch}{0.55}
 \centering
 \caption{Basis function selection for the SW dual}
        \begin{tabular}{c c rrr c rrr c rrr}
    	\toprule
    	 & & \multicolumn{3}{c}{\# Basis functions} & & \multicolumn{3}{c}{Time (s)} & & \multicolumn{3}{c}{Bound (scaled)} \\
	 \cmidrule{3-5}
	 \cmidrule{7-9}
	 \cmidrule{11-13}
      	      Option &  & \multicolumn{1}{c}{$T=4$} & \multicolumn{1}{c}{$T=6$} & \multicolumn{1}{c}{$T=8$}&  & \multicolumn{1}{c}{$T=4$} & \multicolumn{1}{c}{$T=6$} & \multicolumn{1}{c}{$T=8$} & & \multicolumn{1}{c}{$T=4$} & \multicolumn{1}{c}{$T=6$} & \multicolumn{1}{c}{$T=8$} \\
		\cmidrule(l){1-1} \cmidrule(l){3-5}  \cmidrule(l){7-9}  \cmidrule(l){11-13}
    	1     & & 315   & 750   & 1365  & & 353.8 & 1412.5 & 6273.9 & & 60.7  & 57.5  & 56.6 \\
    	2     & & 180   & 300   & 420   & & 91.7  & 243.7 & 632.5 & & 57.9  & 55.4  & 53.8 \\
    	3     & & 135   & 300   & 525   & & 109.0   & 288.7 & 795.9 & & 58.1  & 55.9  & 54.1 \\
    	4     & & 90    & 150   & 210   & & 157.2 & 795.9 & 323.6 & & 54.7  & 54.2  & 46.1 \\
    	\bottomrule
    \end{tabular}%
 \label{tab:bf-options-sw}%
\end{table}%

In Table \ref{tab:bf-options-sw}, Option 1 has the largest number of basis functions and LDDR variables, and the highest lower bound among all the four options. We also carried out a pairwise comparison between Option 1 and the others by performing a $ t $-test. The test confirmed that the difference in the means is statistically significant (at 95\% confidence). The higher quality of the bound using Option 1 comes at the price of a larger solution time. If the computation becomes cumbersome, the next candidates are Options 2 and 3 which do not have a significant difference in their means.

\begin{table}[tbp]
\renewcommand{\arraystretch}{0.55}
 \centering
 \caption{Basis function selection for the NA dual}
     \begin{tabular}{c c rrr c rrr c rrr}
    	\toprule
    	 & & \multicolumn{3}{c}{\# Basis functions} & & \multicolumn{3}{c}{Time (s)} & & \multicolumn{3}{c}{Bound (scaled)} \\
	 \cmidrule{3-5}
	 \cmidrule{7-9}
	 \cmidrule{11-13}
      	      Option &  & \multicolumn{1}{c}{$T=4$} & \multicolumn{1}{c}{$T=6$} & \multicolumn{1}{c}{$T=8$}&  & \multicolumn{1}{c}{$T=4$} & \multicolumn{1}{c}{$T=6$} & \multicolumn{1}{c}{$T=8$} & & \multicolumn{1}{c}{$T=4$} & \multicolumn{1}{c}{$T=6$} & \multicolumn{1}{c}{$T=8$} \\
		\cmidrule(l){1-1} \cmidrule(l){3-5}  \cmidrule(l){7-9}  \cmidrule(l){11-13}
  	1     & & 756   & 1800  & 3276  & & 76.9  & 403.2 & 1800.8 & & 99.8  & 99.9  & 99.8 \\
  	2    &  & 432   & 720   & 1008  & & 44.9  & 232.3 & 770.5 & & 93.4  & 94.2  & 93.1 \\
  	3    &  & 324   & 720   & 1260 &  & 61.1  & 250.0  & 776.8 & & 99.8  & 99.3  & 99.5 \\
  	4    &  & 216   & 360   & 504   & & 34.2  & 219.5 & 735.3 & & 93.3  & 93.6  & 93.1 \\
  	\bottomrule
  \end{tabular}%
 \label{tab:bf-options-na}%
\end{table}%
The results for the NA dual are presented in Table \ref{tab:bf-options-na}. In this case we find that Options 1 and 3, which use all future time periods rather than just the next time period, have clearly better quality of the bound. We conducted a pairwise $ t $-test to determine if the differences are statistically significant. In particular, we compared all the options against Option 3, which seems to return  a very high quality solution in a reasonable amount of time. The test revealed that there is no statistically distinguishable difference between the bounds obtained by Option 1 and Option 3.  Option 3 is statistically significantly better than Options 2 and 4. This shows that having information about the whole planning horizon is beneficial, rather than just considering the next stage ahead, but for the NA dual there was no observed benefit to using demands for products different from the one being relaxed in the NA constraint in the basis functions. \rev{Note that these conclusions depend on the specific instances we are solving and might differ for another set of instances}. In the light of the above discussions, in the rest of the experiments, Options 1 and 3 are used for the SW and NA duals, respectively.

\subsubsection{NA Variable Selection}\label{sec:na-var-selection}
\revtwo{In the NA problem reformulation, another possibility to reduce the restricted dual problem size is to fix the dual variables associated with some of the NA constraints \eqref{eq:const-reform-nonanticipative} to zero.   
Table \ref{tab:na-var-selection} examines this by considering three options \rev{for instances with $J=3, \rho = 0.6, \rho^Y = 0.2$}. Results show that the bounds from penalizing the NA constraints on only the $ x $ variables are not significantly worse than those obtained by penalizing the NA constraints on all variables.  Since the former yields lower bounds that are indistinguishable from the best, in less time, we use this option in the remaining experiments. A more detailed discussion is provided in the Appendix.}

\rev{
\begin{table}[tbp]
\renewcommand{\arraystretch}{0.55}
 \centering
 \caption{NA variable selection}
\begin{tabular}{l rrr rrr}
\toprule
 & \multicolumn{3}{c}{Time (s)} & \multicolumn{3}{c}{Bound} \\
\cmidrule(l){2-4}\cmidrule(l){5-7} Option     & \multicolumn{1}{c}{$T=4$} & \multicolumn{1}{c}{$T=6$} & \multicolumn{1}{c}{$T=8$} & \multicolumn{1}{c}{$T=4$} & \multicolumn{1}{c}{$T=6$} & \multicolumn{1}{c}{$T=8$} \\
\cmidrule(l){1-1}\cmidrule(l){2-4}\cmidrule(l){5-7}
$x$   & 61.1  & 250.0   & 776.8 & 99.8  & 99.3  & 99.5 \\
$i^+, i^-$ & 68.1  & 188.0   & 653.1 & 92.5  & 91.0  & 90.4 \\
$x, i^+, i^-$ & 175.4 & 849.9 & 1887.3 & 100.0 & 100.0 & 100.0 \\
\bottomrule
\end{tabular}%
 \label{tab:na-var-selection}%
\end{table}%
}

\subsubsection{Bound Comparison}
\begin{figure}[htbp]
	\centering
	\begin{subfigure}[t]{0.5\textwidth}
		\centering
		\includegraphics[width=\linewidth]{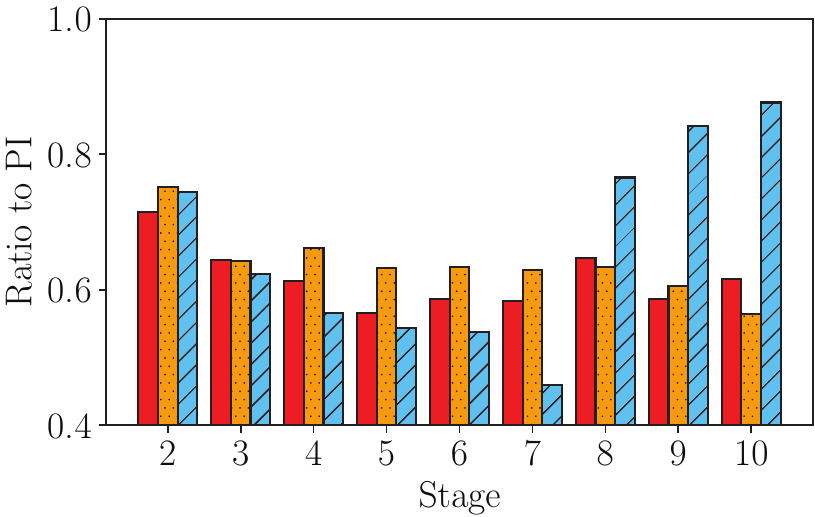}
		\caption{SW lower bound}
		\label{fig:swlb}
	\end{subfigure}
	\hfill
	\begin{subfigure}[t]{0.49\textwidth}
		\centering
		\includegraphics[width=\linewidth]{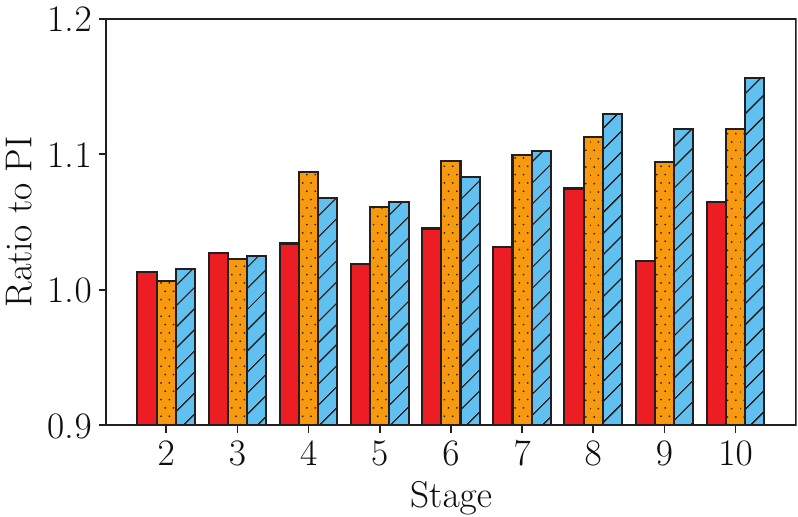}
		\caption{NA lower bound}
		\label{fig:nalb}
	\end{subfigure}
	\caption{Performance of lower bounding algorithms with respect to each other and the PI. Solid bars represent instances with $ J = 3, \rho=0.2, \rho^Y=0.6 $, dotted bars are instances with $ J = 3, \rho=0.6, \rho^Y=0.2 $, and in dashed bars we have $ J = 6, \rho=0.6, \rho^Y=0.2 $.}
	\label{fig:LBs}
\end{figure}

\begin{figure}[btp]
	\centering
	\begin{subfigure}[t]{0.49\textwidth}
		\centering
		\includegraphics[width=\linewidth]{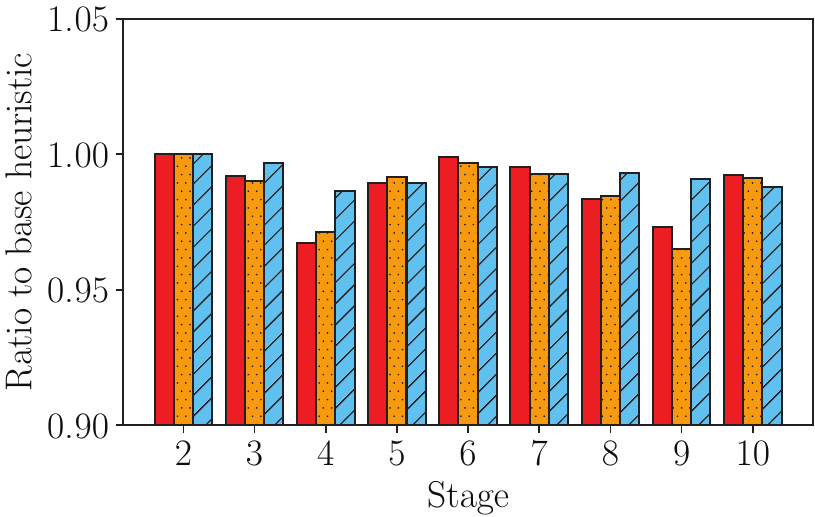}
		\caption{SW upper bound}
		\label{fig:swub}
	\end{subfigure}
	\hfill
	\begin{subfigure}[t]{0.49\textwidth}
		\centering
		\includegraphics[width=\linewidth]{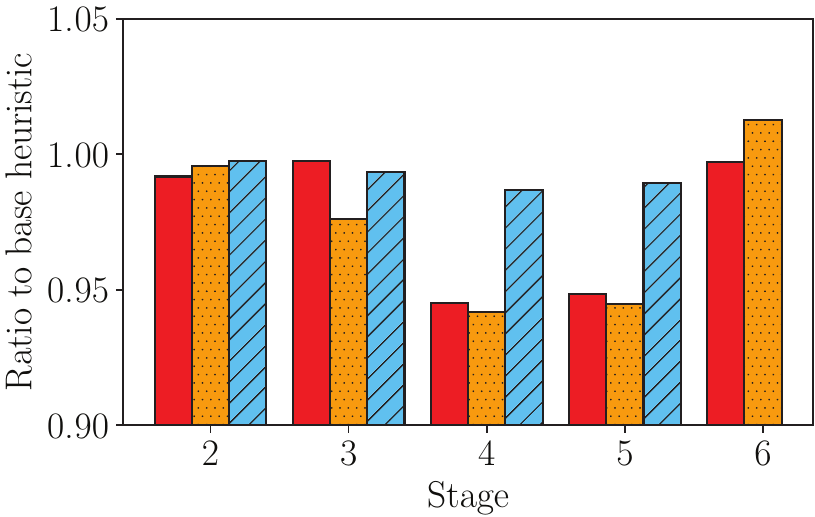}
		\caption{NA upper bound}
		\label{fig:swubj}
	\end{subfigure}
	\caption{Performance of upper bounding algorithms with respect to each other and the conditional expected value policy. Solid bars represent instances with $ J = 3, \rho=0.2, \rho^Y=0.6 $, dotted bars are instances with $ J = 3, \rho=0.6, \rho^Y=0.2 $, and in dashed bars we have $ J = 6, \rho=0.6, \rho^Y=0.2 $.}
	\label{fig:UBs}
\end{figure}

We now compare the quality of the bounds returned by the two Lagrangian duals, and their respective \revtwo{primal} bounds, over a variety of test instances. Figures \ref{fig:LBs} and \ref{fig:UBs} illustrate the performance of the lower \revtwo{(dual)} and upper \revtwo{(primal)} bounding algorithms \revtwo{(see Section \ref{sec:implementation}),} for instances with different stages, product types, and levels of stochasticity in demands. The reported numbers are normalized as follows. The lower bounds obtained for each instance are divided by the PI lower bound for that instance, so that the reported value is the lower bound relative to the PI lower bound and values over 1.0 indicate an improvement in the lower bound. The upper bounds are divided by the bound returned by the conditional expected value policy, so that a value below 1.0 indicates an improvement in the upper bound. For both cases the actual numbers are available in
the Appendix.

As shown by Figure \ref{fig:nalb}, the NA dual has been able to improve upon the PI lower bound in all instances. This improvement generally becomes more evident as the number of stages grows. From Figure \ref{fig:nalb}, it can be seen that the NA dual is performing better with a higher level of stochasticity (less correlation in the autoregressive process and a higher weight for variation). This could be due to a stronger PI bound when the system is less volatile. On the other hand, as seen in Figure \ref{fig:swlb}, the SW dual provides significantly lower bounds (as expected from Theorem \ref{thm:sw-na-compare}) than those from the NA dual, and even is significantly worse than the PI bound on these test instances. 

Figure \ref{fig:UBs}  presents the results comparing the upper bounds obtained from the dual-driven policies to that from the conditional expected value policy. We see from Figure \ref{fig:swub} that the SW dual-driven policy provides modest but consistent improvement over the conditional expected value policy. For the NA dual-driven policy we report results only for stages $T \leq 6$, since estimating the value of the policy was too time-consuming for larger instances
\rev{ and takes longer than our maximum running time of 6 hours. It should be noted that this is caused by the very large number of scenarios in the evaluation sample, and on average obtaining solutions for one sample path does not take more than 45 seconds in the largest instances. The solution times of our lower and upper bounding techniques can be found in the Appendix.} We find that for instances with $J=3$ and $3 \leq T \leq 5$ the NA dual-driven policy was able to significantly improve upon the conditional expected value policy. We suspect that for instances with larger $T$ or $J$, the lack of improvement is due to a poor approximation of the two-stage stochastic MIP \eqref{eq:twostage-master} that is solved in this method, as we approximated this problem with just 25 scenarios considering the computational limitations. Although we do not explore this here, the NA dual-driven policy may be more practical if a decomposition algorithm is used to solve the problem \eqref{eq:twostage-master}, enabling the use of many more scenarios to approximate it.
\begin{figure}[htbp]
	\centering
	\begin{subfigure}[t]{0.49\textwidth}
		\centering
		\includegraphics[width=\linewidth]{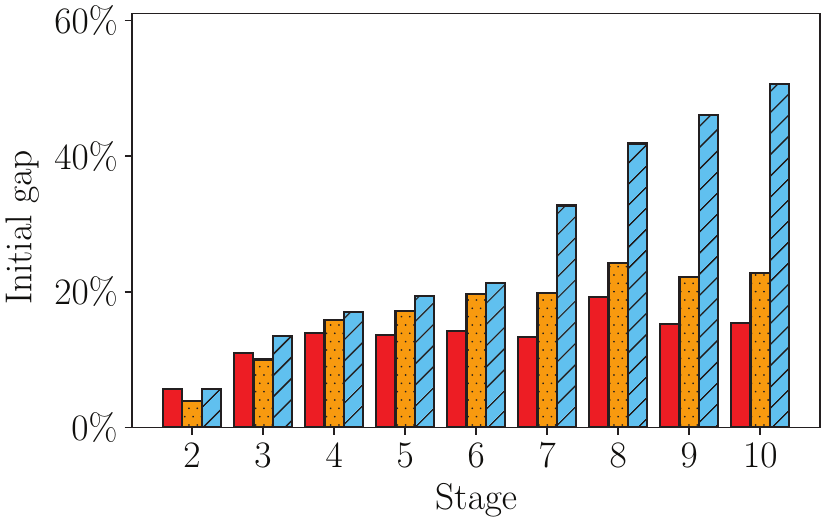}
		\caption{Gap between PI lower bound and conditional expected value upper bound for all instances}
		\label{fig:initgap}
	\end{subfigure}
	\hfill
	\begin{subfigure}[t]{0.5\textwidth}
		\centering
		\includegraphics[width=\linewidth]{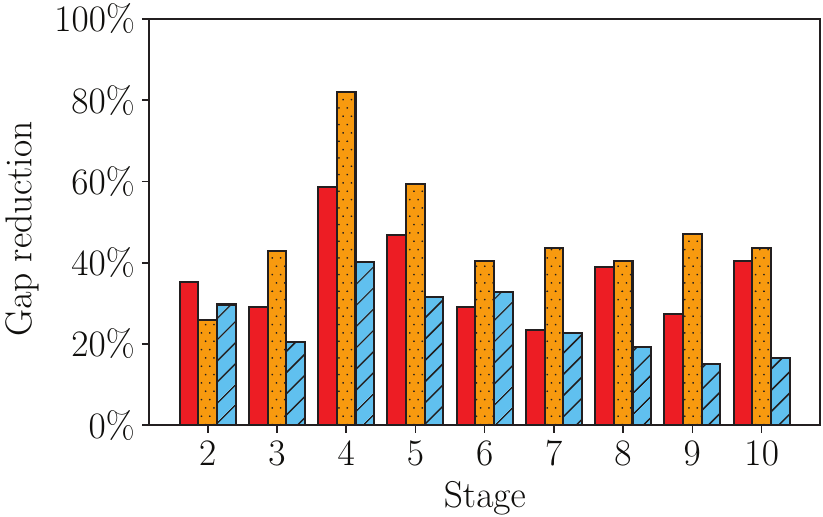}
		\caption{Gap reduction obtained from combining the improvements of lower and upper bounding techniques}
		\label{fig:gap}
	\end{subfigure}
	\caption{Gap reduction percentage with respect to the gap between the PI and the conditional expected value policy. Solid bars represent instances with $ J = 3, \rho=0.2, \rho^Y=0.6 $, dotted bars are instances with $ J = 3, \rho=0.6, \rho^Y=0.2 $, and in dashed bars we have $ J = 6, \rho=0.6, \rho^Y=0.2 $.}
	\label{fig:gaps}
\end{figure}

Figure \ref{fig:initgap} presents the gap between the PI lower bound and the conditional expected value upper bound across various instances, which ranges from 3.91\% to 50.67\%. In Figure \ref{fig:gap} we present the fraction of this gap that is closed using the combined improvements from the NA dual lower bound and the SW dual upper bound. We find that the gap is reduced over all instances, and is more pronounced on instances with $J=3$ and instances with $\rho=0.6,\rho^Y=0.2$ (the instances in which the demands have higher variability).
\rev{
\subsection{Stochastic Routing and Wavelength Assignment Problem}\label{sec:MSIP-RWA}
In a common optical telecommunications network, links are optical fibers that can simultaneously transmit signals on different wavelengths. The routing and wavelength assignment (RWA) problem is the classical resource allocation problem for such networks. 
Given a network, a set of wavelengths, and a set of connection requests (origin and destination node pairs), the RWA problem assigns a path and a wavelength to each request, with the objective of maximizing the number of granted requests. This pair of path and wavelink is called a \emph{lightpath}. The building block of lightpaths is a \emph{wavelink}, the pair of a link and a wavelength. Wavelinks on a valid lightpath all have the same wavelength (\emph{wavelength continuity}). To avoid \emph{wavelength conflict}, no two lightpaths can share a wavelink.
Figure \ref{fig:wdm-network} presents a bidirectional optical network with three wavelengths (represented by different colors). While the dotted red and wavy blue lightpaths are valid, the green dashed lightpaths are invalid as they share the wavelink from node 1 to 3. 
\begin{figure}[h]
\centering
\input{./Figures/Abilene.tex}%
\caption{\rev{An optical network with valid and invalid lightpaths}}
\label{fig:wdm-network}
\end{figure}  

Due to the scale of the problems in telecommunications, most solution methods are heuristics, an overview of which is compiled by \cite{simmons2020routing}, with scalable exact algorithms only recently proposed by \cite{jaumard2017efficient} for deterministic traffic. For stochastic incoming traffic, \cite{daryalal2020stochastic} proposed a two-stage stochastic MIP model and solved it via a decomposition framework. In reality, the set of connection requests arrives randomly over time, suggesting that a multistage model might better capture the dynamic process.

Consider an optical network with a given physical topology represented by a multigraph $\Graph = (\nodeSet,\edgeSet) $ and a set of available wavelengths $\waveSet$ indexed by $\wavelength$. We denote by $\edgeSetw\subseteq\edgeSet$ the set of free fiber links on wavelength $\wavelength$. We assume requests remain on the network throughout the planning horizon, and let $ \DnewMSIP $ denote a random parameter representing the number of requested lightpaths between node pair $(s,d)$ at stage $ t $. At stage $t$, we let
$\linksceMSIP$ be a binary (recourse) decision variable equal to 1 if wavelink $ (\edge,\wavelength) $  is used for a new connection request between $ (s,d)\in \SDsceMSIP=\{(s,d) \in \nodeSet\times \nodeSet: \DnewMSIP > 0\}$ and equal to 0 otherwise, and let $ \wavelinksceMSIP $ be a binary (state) variable equal to 1 if wavelink $ (\edge,\wavelength) $  is serving any node pair at stage $ t $ (based on a current or previous accepted connection). Finally,
$\numgrantedsceMSIP$ is a (recourse) decision variable measuring the number of granted requests between node pair $(s,d)$ during stage $t$. With these decision variables, the \MSMIP{} model for multistage  RWA (MSRWA) is:
\BSE
\label{eq: SmaxRWA_MSIP}
\begin{alignat}{5}
	\max \ \ & \Exp_{\xiT}  \left[\sum_{t\in[T]}\sum_{(s,d) \in \SDsceMSIP} \numgrantedsceMSIP  \right] \label{eq: OF_MSIP}\\
	\text{s.t.}\ \ & 
	\sum_{(s,d) \in \SDcurr} \linkMSIP  = \wavelinkMSIP &\quad& \wavelength\in \waveSet, \edge \in \edgeSetw \label{eq: wavelink_usage_t1}\\
	&\wavelinksceMSIP  = \wavelinkprevsceMSIP + \sum_{(s,d) \in \SDsceMSIP} \linksceMSIP &\ & t\in[2,T], \wavelength\in \waveSet, \edge \in \edgeSetw, \allxi \label{eq: conflict_link_demands_scenarios_MSIP}\\
	&\sum_{\stackrel{\edge \in}{\nodefunc^+(\node)\cap \edgeSetw}} \linksceMSIP = \sum_{\stackrel{\edge \in}{\nodefunc^-(\node)\cap \edgeSetw}} \linksceMSIP && t\in[T], \wavelength \in \waveSet, (s,d) \in \SDsceMSIP,  \label{eq: flow_conservation_scenarios_MSIP} \\*[-5mm]
	&&& \node\in \nodeSet\setminus\{s,d\}, \allxi \nonumber \\
	& \sum_{\wavelength \in \waveSet} \sum_{\stackrel{\edge \in}{\nodefunc^-(s)\cap \edgeSetw}}\linksceMSIP  =0 && t\in[T], (s,d) \in \SDsceMSIP, \allxi\label{eq: loop_elimination1_scenarios_MSIP}\\
	& \sum_{\wavelength \in \waveSet} \sum_{\stackrel{\edge \in}{\nodefunc^+(d)\cap \edgeSetw}}  \linksceMSIP = 0 && t\in[T], (s,d) \in \SDsceMSIP, \allxi\label{eq: loop_elimination2_scenarios_MSIP}\\
	&\numgrantedsceMSIP = \sum_{\wavelength \in \waveSet}\sum_{\stackrel{\edge \in}{\nodefunc^+(s)\cap \edgeSetw}} \linksceMSIP && t\in[T], (s,d) \in \SDsceMSIP,\allxi\label{eq: granted_scenarios_MSIP}\\
	&\numgrantedsceMSIP \leq \DnewMSIP &&t\in[T], (s,d) \in \SDsceMSIP, \allxi \label{eq: demand_UB_scenarios_MSIP} \\
	& \linksceMSIP \in\{0,1\} && t\in[T], (s,d) \in \SDsceMSIP,\wavelength\in\waveSet, \edge\in\edgeSetw,\label{eq: linksce_bounds_MSIP}\\
	&&& \allxi \nonumber\\
	& \wavelinksceMSIP \in\{0,1\} && t\in[T],\wavelength\in\waveSet, \edge\in\edgeSetw,\allxi\label{eq: wavelinksce_bounds_MSIP}
\end{alignat}
\ESE
where $\nodefunc^{+(-)}\subseteq\edgeSet$ is the set of outgoing (incoming) links of node $\node\in\nodeSet$. The objective in \eqref{eq: OF_MSIP} maximizes the expected number of granted  requests. 
Constraints \eqref{eq: wavelink_usage_t1} and \eqref{eq: conflict_link_demands_scenarios_MSIP} are state equations that keep track of whether or not each wavelink has been used by stage $t$, which together with the binary restrictions on the state variables	 prevent a wavelink from being used more than once in the planning horizon.
Constraints \eqref{eq: flow_conservation_scenarios_MSIP} -\eqref{eq: granted_scenarios_MSIP} ensure that a valid lightpath is created for each accepted connection. 
Constraints \eqref{eq: demand_UB_scenarios_MSIP} bound the number of lightpath assignments to a node pair by its demand.
Constraints \eqref{eq: linksce_bounds_MSIP}-\eqref{eq: wavelinksce_bounds_MSIP} are the integrality and bound restrictions.

\subsubsection{Data and Basis Functions} \label{sec:data-bf-RWA}
For the underlying topology of our instances, we have considered \abilene{}, a standard long-haul network in the literature of the RWA problem \citep{orlowski2010sndlib}, with $|\waveSet| \in \{1,3\}$. The total number of request arrivals in the incoming traffic $\traffic_t$ follows an integer-valued GARCH process (INGARCH, \cite{ferland2006integer}), where the deviates have a Poisson distribution with a parameter $\lambda_t$ \citep{chen2015spectrum,xiong2018sdn} that depends on $\traffic_{t-1}$ and $\lambda_{t-1}$:
$$ \left\lbrace
\begin{array}{lr}
\traffic_t | \mathcal{F}_{t-1} : \text{Poisson}(\lambda_t)\\
\lambda_t = \gamma_0 + \gamma_1 \traffic_{t-1} + \delta_1\lambda_{t-1}, & \hspace{1cm}  \gamma_0 > 0, \gamma_1\geq 0, \delta_1\geq 0, \gamma_1+\delta_1 <1
\end{array}
\right.
 $$
and the traffic is distributed uniformly among the node pairs \citep{giorgetti2015dynamic} to yield the random variables $\DnewMSIP$. We can show that the expectation of the process is $\bar{\traffic} = \displaystyle\frac{\gamma_0}{1-(\gamma_1 + \delta_1)}$ and the $k$-step ahead conditional expectation is
$$\Exp[\traffic_{t+1}|\mathcal{F}_{t}\ ] = \gamma_0 \sum_{i=0}^{k-2}(\gamma_1 + \delta_1)^i+(\gamma_1 + \delta_1)^{k-1}(\gamma_0 + \gamma_1 \traffic_{t} + \delta_1\lambda_{t}).$$
For $T \in\{3,4,5\}$ we have generated instances with  $\bar{\traffic}\in\{10,30\}$, $\gamma_1=0.5$ and two levels of correlation between $\lambda_t$ parameters of consecutive stages $\delta_1\in\{0.2,0.05\}$.  
\revtwo{For the NA variable selection we have two options: state variables $\wavelinksceMSIP$ and recourse variables $\linksceMSIP$. Since by using $\linksceMSIP$ the number of LDDR decision variables would grow with the number of connection requests, we opt for  $\wavelinksceMSIP$ variables.  Table \ref{tab:options-rwa} presents our considered options for the basis functions. In Option 2, we try to reduce the number of basis functions by selecting the node pairs that have a higher chance of being served by a wavelink. For a fixed wavelink $(\hat{\edge} , \hat{\wavelength}) $, we consider the demands of the node pairs that belong to the  set
$ \SDlsceMSIPBF = \Big\{ (s,d) \in \nodeSet\times\nodeSet : \DnewMSIPBF > 0,\ \hat{\edge} \in \bigcup_{i=1}^2\textsc{SP}^i(s,d) \Big\}, $ where $ \textsc{SP}^i(s,d)  $ is the set of $ i^{\text{th}} $ shortest paths between $s$ and $d$. 
}
\begin{table}[h]
\renewcommand{\arraystretch}{0.55}
	\centering
	\caption{\rev{Options to be used in basis functions, at stage $ \hat{t} $ for wavelink $(\hat{\edge} , \hat{\wavelength}) $}}
	\begin{tabular}{c cll l cll l}
		\toprule
		Option & & \multicolumn{3}{c}{SW Dual} & & \multicolumn{3}{c}{NA dual} \\
		\midrule
		1     & & $1 \ ; \displaystyle\sum_{(s,d)\in \SDsceMSIPBF}\DnewMSIPBF$, & & $\forall t\leq \hat{t}$ &&  $1\ ; \DnewMSIPBF$, & $(s,d)\in \SDsceMSIPBF$, & $\forall t\leq T: \hat{t} < t$ \\
		2     & & $1 \ ; \DnewMSIPBF$, & $(s,d)\in \SDlsceMSIPBF$, & $\forall t\leq \hat{t}$ &&  $1\ ; \DnewMSIPBF$, & $(s,d)\in \SDlsceMSIPBF$, & $\forall t\leq T: \hat{t} < t$ \\
		\bottomrule
	\end{tabular}%
	\label{tab:options-rwa}%
\end{table}%

\subsubsection{Bound Comparison}\label{sec:bounds-RWA}
In this section, we present the results of our \revtwo{primal and dual} bounding methods \revtwo{(to obtain lower and upper bounds for the maximization problem, respectively)} on the MSRWA instances using basis function Options 1 and 2 (Table \ref{tab:options-rwa}). Since our best results in the previous section were obtained by using the SW dual to drive a primal policy, and the NA dual to obtain lower bounds, we only consider these options for the RWA problem. 

The results are presented in Table \ref{tab:bounds-rwa}. \revtwo{Note that, contrasted with the MSLot problem, the MSRWA is a maximization problem, thus the primal and dual bounding techniques of Section \ref{sec:implementation} lead to lower and upper bounds, respectively.} For each option, in the first column we report the ratio of the NA upper bound to the PI upper bound, in the second column we report the ratio of the SW lower bound to the conditional expected value lower bounding policy, and in the third column we report the overall gap reduction (relative to the initial gap) obtained by considering both bounds.  Using Option 1, the initial gap is considerably improved, with gap reductions varying from 7\% ($10^{\text{th}}$ instance) to 63\% ($1^{\text{st}}$ instance). Our upper and lower bounding methods both contribute to this gap reduction. Although the upper bounds from Option 1 should theoretically always be better than those from Option 2 (because the set of basis functions in Option 2 is a subset of those in Option 1), the higher difficulty in solving the Lagrangian dual for Option 1 meant that we needed to use a smaller sample size, often resulting in worse bounds than were obtained with Option 2 when evaluated out of sample. \revtwo{The solution times of the bounding techniques for each basis function option are presented in the Appendix.}
As evidenced by the results in Table \ref{tab:bounds-rwa}, using Option 2, even with a subset of the standard basis functions we are able to reduce the gaps up to 15\% more in comparison to the Option 1. This demonstrates the importance of the design of basis functions. \revtwo{Another important observation is the impact of the temporal dependency present in the instances on the quality of the obtained bounds. More gap is reduced in instances with larger values of $\delta_1$ (representing the correlation level between consecutive stages), and both lower and upper bounding techniques demonstrate better performance for instances with higher correlation among stages. We also found that the initial gap and gap reduction tend to be higher for instances with more stages $T$.}
\begin{table}[htbp]
\renewcommand{\arraystretch}{0.55}
  \centering
  \caption{\rev{Gap reductions for MSRWA instances, obtained from combining the improvements of NA upper bound and SW lower bound.}}
 \small{
\begin{tabular}{lc c cc c cc c}
\toprule
      &       &       & \multicolumn{3}{c}{Basis Function Option 1} & \multicolumn{3}{c}{Basis Function Option 2} \\
\cmidrule(lr){4-6} \cmidrule(lr){7-9} $(|\waveSet|,\bar{\traffic},\gamma_0,\delta_1)$ & $T$   & \multicolumn{1}{c}{Initial Gap} & \multicolumn{1}{c}{$\frac{\text{NA UB}}{\text{PI}}$} & \multicolumn{1}{c}{$\frac{\text{SW LB}}{\text{Exp}}$} & \multicolumn{1}{c}{Gap reduction} & \multicolumn{1}{c}{$\frac{\text{NA UB}}{\text{PI}}$} & \multicolumn{1}{c}{$\frac{\text{SW LB}}{\text{Exp}}$} & \multicolumn{1}{c}{Gap reduction} \\
\cmidrule(lr){1-1} \cmidrule(lr){2-2} \cmidrule(lr){3-3}\cmidrule(lr){4-6} \cmidrule(lr){7-9} 
\multirow{3}[-2]{*}{$(1,10,3,0.2)$} & 3     & 27.0\% & 0.94  & 1.08  & 63.0\% & 0.94  & 1.05  & 50.4\% \\
      & 4     & 43.1\% & 0.95  & 1.09  & 41.8\% & 0.93  & 1.12  & 56.3\% \\
      & 5     & 45.1\% & 0.97  & 1.15  & 50.1\% & 0.96  & 1.16  & 54.5\% \\
\midrule
\multirow{3}[-2]{*}{$(1,10,4.5,0.05)$} & 3     & 22.0\% & 0.97  & 1.02  & 27.4\% & 0.97  & 1.05  & 42.8\% \\
      & 4     & 34.8\% & 0.97  & 1.06  & 32.0\% & 0.95  & 1.07  & 41.8\% \\
      & 5     & 51.4\% & 0.97  & 1.18  & 51.4\% & 0.94  & 1.18  & 59.8\% \\
\midrule
\multirow{3}[-2]{*}{$(3,30,9,0.2)$} & 3     & 27.7\% & 0.99  & 1.06  & 30.5\% & 0.99  & 1.08  & 39.9\% \\
      & 4     & 24.2\% & 0.98  & 1.03  & 22.5\% & 0.98  & 1.04  & 29.1\% \\
      & 5     & 34.0\% & 0.98  & 1.09  & 42.2\% & 0.98  & 1.10  & 44.3\% \\
\midrule
\multirow{3}[-2]{*}{$(3,30,13.5,0.05)$} & 3     & 20.3\% & 1.00  & 1.01  & \hspace*{0.15cm}7.0\% & 1.00  & 1.03  & 18.0\% \\
      & 4     & 18.3\% & 0.99  & 1.03  & 27.4\% & 0.99  & 1.03  & 29.0\% \\
      & 5     & 32.2\% & 0.99  & 1.03  & 17.7\% & 0.98  & 1.02  & 15.6\% \\
\bottomrule
\end{tabular}%
}
  \label{tab:bounds-rwa}%
\end{table}%
\spacing
}

\section{Conclusion} 
\label{sec:conclusion} 
In this work, we introduced the idea of Lagrangian dual decision rules where decision rules are used in the Lagrangian dual of an MSMIP. The result is an approximation problem that can be solved by stochastic approximation or sample average approximation. The approximate problem does not have a multistage structure, and hence does not require a scenario tree for its approximate solution. Two lower bounding policies based on two Lagrangian duals are proposed: stagewise and nonanticipative, where the former is an easier problem consisting of single period subproblems, while the latter can potentially lead to better bounds if the basis functions are selected properly.  The solutions to both of these duals can be incorporated in constructing primal policies. The lower and upper bounding methods were evaluated by solving instances of \rev{two problem classes:} a multi-item stochastic lot-sizing problem \rev{and routing and wavelength assignment problem in optical networks}. The results show that our methods can substantially reduce the optimality gap relative to the use of two general-purpose bounding policies. Future work includes design of more scalable dual-driven policies and investigation of other Lagrangian dual decision rule structures, such as piecewise linear form.

\ACKNOWLEDGMENT{This work was supported by Natural Sciences and Engineering Research Council of Canada [Grant RGPIN-2018-04984], NSF award CMMI-1634597 and by the Department of Energy, Office of Science, Office of Advanced Scientific Computing Research, Applied Mathematics program under Contract Number DE-AC02-06CH113.}

\bibliography{Bibilio/literature.bib}
\bibliographystyle{informs2014}

\begin{appendices}
	
In the following, we have provided further details on the proofs and results given in the paper.
\section{Proof of Lemma \ref{lem:na-reform}} \label{sec:proof-lem-1}
Expanding the objective function in \eqref{eq:obj-lagrange-na-t} and rearranging the terms yields
	\begin{equation}\label{eq:expanded_expectation_NA}
	\sum_{t\in \set{T}} \left( \Exp_{\xiT} \left[c_t(\xisupt)^\top y_t(\xiT) \right] + \Exp_{\xiT} \left[ \gamma_t(\xiT)^\top y_t(\xiT) \right] - \Exp_{\xiT} \left[ \gamma_t(\xiT)^\top \Exp_{\xiprimeT}[y_t(\xiprimeT) | \xiprimesupt = \xisupt]\right]\right).  
	\end{equation}	
For the last term in \eqref{eq:expanded_expectation_NA}, the following equalities hold
\bsub
\begin{eqnarray}
\Exp_{\xiT} \left[ \gamma_t(\xiT)^\top \Exp_{\xiprimeT}[y_t(\xiprimeT) | \xiprimesupt = \xisupt]\right] & = & \Exp_{\xiT} \left[  \Exp_{\xiprimeT}[\gamma_t(\xiT)^\top y_t(\xiprimeT) | \xiprimesupt = \xisupt]\right] \label{eq:lem-na-a} \\
& = & \Exp_{\xiprimeT} \left[  \Exp_{\xiT}[\gamma_t(\xiT)^\top y_t(\xiprimeT) | \xiprimesupt = \xisupt]\right]  \label{eq:lem-na-b} \\
& = & \Exp_{\xiprimeT} \left[  \Exp_{\xiT}[\gamma_t(\xiT) | \xiprimesupt = \xisupt]\>y_t(\xiprimeT)\right]  \label{eq:lem-na-c}\\
& = & \Exp_{\xiT} \left[ \Exp_{\xiprimeT}[\gamma_t(\xiprimeT) | \xiprimesupt = \xisupt] y_t(\xiT) \right]. \label{eq:lem-na-d}
\end{eqnarray}
\esub
First note that, in $ \Exp_{\xiT} \left[ \gamma_t(\xiT)^\top \Exp_{\xiprimeT}[y_t(\xiprimeT) | \xiprimesupt = \xisupt]\right] $, since the dual function $ \gamma_t(\xiT) $ is fixed, $ \gamma_t(\xiT) $ is a vector of numbers inside the first expectation, hence it can be pushed inside the second expectation (Equation \eqref{eq:lem-na-a}).
As we have assumed that $ \gamma_t(\xiT) $ is a member of the set $ \Gamma_t $ with its expectation being bounded, and \revtwo{$\Exp[\dia{Y_t(\xiprimesupt)}]$, the expectation of the} diameter of the set $ Y_t(\xiprimesupt) $ to which $ y_t(\xiprimesupt) $ belongs, is finite, using Fubini-Tonelli Theorem \citep{knapp2005basic} the order of the two expectations can be exchanged (Equation \eqref{eq:lem-na-b}).
Then we can take $ y_t(\xiprimeT) $ out from the inside expectation, as it is just a vector of numbers inside (Equation \eqref{eq:lem-na-c}). Since $ \xiT $ and $ \xiprimeT $ have the same support and distribution, the last equality (Equation \eqref{eq:lem-na-d}) is satisfied. Equality \eqref{eqq:lemma1} can then be proven by substitution.

\section{Proof of Lemma \ref{lem:primal-char}} \label{sec:proof-lem-2}
First observe that since the objective in \eqref{eq:zld}  is linear, the restricted Lagrangian dual problem \eqref{eq:rld} can be written as
\begin{align*}
z^{\text{RLD}} = \max_{\lambda, \alpha} \ \min_x \ \ & \{ c^\top x + \lambda^\top(d-Dx) : x \in \conv{X} \}  \\
\text{s.t.} \ \ & \lambda - G\alpha = 0
\end{align*}
Let $ \{x^i\}_{i \in \set{M}} $ and $ \{r^k\}_{k\in\set{K}} $ be the complete set of extreme points and extreme rays of $ \conv{X}$, respectively. Then, for any fixed $ \lambda $, we have
\begin{equation*}
z(\lambda) = \left\lbrace
\begin{array}{ll}
-\infty, & \exists r^k: (c^\top-\lambda^\top D)r^k < 0\\
\displaystyle \min_{i\in\set{M}} \{c^\top x^i + \lambda^\top (d-Dx^i)\}, & \text{otherwise}.
\end{array}
\right.
\end{equation*}
Therefore we can reformulate the Lagrangian dual problem as 
\begin{align*}
z^{\text{LD}} = \max_{\lambda, \alpha} \ \min_{i\in\set{M}} \ \ & \{c^\top x^i + \lambda^\top (d-Dx^i)\} \\
\text{s.t.} \ \ & \lambda - G\alpha = 0 \\
& (c^\top -\lambda^\top D)r^k \geq  0, \quad k\in\set{K}    
\end{align*}
which is equivalent to
\bsub
\label{eq:lem-proof}
\begin{align}
z^{\text{LD}} = \max_{\lambda, \alpha, \eta} \ \ &\eta  \\
\text{s.t.} \ \ & { \lambda - G\alpha = 0} && \hspace*{2.5cm} && { (\theta)} \label{eq:lem-proof-affine} \\
& \lambda^\top Dr^k \leq  c^\top r^k && \quad k\in\set{K} && (\beta_k)\label{eq:lem-proof-const1} \\
& \eta + \lambda^\top (Dx^i - d) \leq c^\top x^i && \quad i\in\set{M} &&  (\gamma_i) \label{eq:lem-proof-const2}
\end{align}
\esub
where $ \theta $, $ \beta_k $ and $ \gamma_i $ are the dual variables   associated with constraints \eqref{eq:lem-proof-affine}, \eqref{eq:lem-proof-const1} and \eqref{eq:lem-proof-const2} respectively.
Now take the dual of the above problem
\begin{align*}
\min_{\theta, \beta, \gamma} \ \ &c^\top \left( \sum_{i\in\set{M}}\gamma_i x^i + \sum_{k\in\set{K}}\beta_k r^k \right)   \\
\text{s.t.} \ \ & \sum_{i\in\set{M}} \gamma_i = 1 \\
& D\left(\sum_{i\in\set{M}}\gamma_i x^I + \sum_{k\in\set{K}}\beta_k r^k \right) + { \theta} = d \\
& { -G^\top  \theta = 0} \\ 
& \gamma, \beta \geq 0
\end{align*}
As we know that $ \conv{X} = \displaystyle \left\{ \sum_{i\in\set{M}}\gamma_i x^i + \sum_{k\in\set{K}}\beta_ kr^k : \sum_{i\in\set{M}} \gamma_i = 1 , \gamma_i, \beta_k \geq 0, \  i\in\set{M}, k \in \set{K} \right\} $, we have
\begin{align}
\min_{ \revtwo{x,\theta}} \ \  & c^\top x \nonumber   \\
\text{s.t.} \ \ & x\in\conv{X} \nonumber \\
& Dx + { \theta} = d  \label{eq:thetadef} \\
& {-G^\top  \theta = 0}. \nonumber
\end{align}
Eliminating the $\theta$ variables using \eqref{eq:thetadef} yields the result.

Note that if $\theta = \boldsymbol{0}$, then the original MIP solution, say $x^*$, is feasible. So, this gives a relaxation of the original MIP.

\section{The Lot-sizing Problem}


\subsection{Stagewise Lagrangian Dual}
Relax state equations \eqref{eq:lot-sizing-state-equations}, except for $ t=1 $. Then, we have
\bsub
	\begin{align}
	\min  \ \ & \Exp
	\left[
	L(\lambda) 
	\right]
	\label{eq:obj-lagrange-lot-sizing-lag}  \\
	\text{s.t.} \ \ & i^-_{1j}(\xiparrand^1)-i^+_{1j}(\xiparrand^1) = d_{1j}(\xiparrand^1),\quad j \in \set{J}, \allxi \\
	& \eqref{eq:lot-sizing-overtime}-\eqref{eq:lot-sizing-integrality}
	\end{align}
\esub
where
\begin{alignat*}{3}
L(\lambda) 
:=  
& \quad 
\sum\limits_{t \in \set{T}} \bigg( \sum\limits_{j \in \set{J}} \Big(
C_{tj}^{i^+}(\xisupt)i^+_{tj}(\xisupt) + C_{tj}^{i^-}(\xisupt)i^-_{tj}(\xisupt) + C_{tj}^{y}(\xisupt)y_{tj}(\xisupt)
\Big) + C_{t}^{o}(\xisupt)o_{t}(\xisupt)
\bigg)\\
& \quad \qquad + 
\sum\limits_{t \in \set{T}} \sum\limits_{j \in \set{J}}\lambda_{tj}(\xisupt)\left(i^-_{tj}(\xisupt)-i^+_{tj}(\xisupt) + i^+_{t-1,j}(\xisuptminusone) - i^-_{t-1,j}(\xisuptminusone) + x_{t-1,j}(\xisuptminusone) - D_{tj}(\xisupt) \right).
\end{alignat*}
Restrict $\lambda_{t}(\xisupt)$ to follow LDDR
\begin{equation*}
\lambda_{tj}(\xisupt) =  \sum\limits_{k \in [\Kt] } \Phi_{tjk}(\xisupt) \alpha_{tjk}.
\end{equation*}
Then, for fixed $\hat{\alpha}$, the objective function \eqref{eq:obj-lagrange-lot-sizing-lag} is equivalent to
\begin{equation}
\sum\limits_{t \in \set{T}}
\Exp
\left[
L_{t}(\hat{\alpha})
\right]
\label{eq:obj-LDR-lot-sizing-lag} 
\end{equation}
where
\begin{alignat*}{4}
& L_{t}(\hat{\alpha}, \xisupt)
:=  \\
& \quad 
\sum\limits_{j \in \set{J}} \Bigg[\Bigg( \sum\limits_{k \in [\Kt]} \Phi_{tjk}(\xisupt)\alpha_{tjk} \Bigg) \big( - D_{tj}(\xisupt) \big)\\
& \quad
+
\Bigg(\Exp \bigg[\sum\limits_{k \in \set{K_{t+1}}} \Phi_{t+1,j,k}(\xisuptplusone)\alpha_{t+1,jk} \> \bigg| \> \xisupt \bigg] \Bigg)x_{tj} 
\\ 
& \quad
+
\Bigg(C_{tj}^{i^+}(\xisupt) - \sum\limits_{k \in [\Kt]} \Phi_{tjk}(\xisupt)\alpha_{tjk} + \Exp \bigg[\sum\limits_{k \in \set{K_{t+1}}} \Phi_{t+1,j,k}(\xisuptplusone)\alpha_{t+1,jk} \> \bigg| \> \xisupt \bigg] \Bigg)i^+_{tj} 
\\ 
& \quad
+
\Bigg(C_{tj}^{i^-}(\xisupt) +\sum\limits_{k \in \set{\Kt}} \Phi_{tjk}(\xisupt)\alpha_{tjk} - \Exp \bigg[\sum\limits_{k \in \set{K_{t+1}}} \Phi_{t+1,j,k}(\xisuptplusone)\alpha_{t+1,jk} \> \bigg| \> \xisupt \bigg] \Bigg)i^-_{tj} 
\\ 
& \quad  
+ C_{tj}^{y}(\xisupt)y_{tj}  \Bigg] + C_{t}^{o}(\xisupt)o_{t}.
\end{alignat*}
For fixed $\hat{\alpha}$, $\xisupt$  and $ t > 1 $, define $\mathcal{L}_t(\hat{\alpha},\xisupt)$ as follows
\bsub
	\label{subeqs:LagDual-lot-sizing-lag}
	\begin{alignat*}{2}
	\mathcal{L}_{t}(\hat{\alpha},\xisupt) \> := \> \min \ \ & 
	L_{t}(\hat{\alpha}, \xisupt)  \\
	\text{s.t.} \ \ 	
	&  \sum\limits_{j \in \set{J}}\big(TS_jy_{tj} +TB_jx_{tj}\big) - o_{t} \leq C_t &\quad&    \\
	& M_{tj}\ y_{tj} - x_{tj} \geq 0, && j \in \set{J}  \\ 
	& i^+_{tj} \leq I_{tj}, &&  j \in \set{J}  \\
	& i^+_{tj} + x_{tj} \leq I_{t+1,j}, &&  j \in \set{J}   \\
	& 0 \leq o_{t} \leq O_{t} &&  \\
	& x_{tj} ,i^+_{tj}, i^-_{tj} \geq 0, && j \in \set{J} \\
	& y_{tj} \in\{0,1\}, && j \in \set{J}. 
	\end{alignat*}
\esub
For $ t = 1 $ we have
\bsub
	\label{subeqs:LagDual-lot-sizing-lag-t1}
	\begin{alignat*}{2}
	\mathcal{L}_{1}(\hat{\alpha}, \xiparrand^1) \> := \> \min \ \ & 
	L_{1}(\hat{\alpha}, \xiparrand^1) \\
	\text{s.t.} \ \ 	
	& i^-_{1j}-i^+_{1j} = d_{1j},&& j \in \set{J}\\
	&  \sum\limits_{j \in \set{J}}\big(TS_jy_{1j} +TB_jx_{1j}\big) - o_{1} \leq C_1 &\quad&    \\
	& M_{1j}\ y_{1j} - x_{1j} \geq 0, && j \in \set{J}  \\ 
	& i^+_{1j} \leq I_{1j}, &&  j \in \set{J} \\
	& i^+_{1j} + x_{1j} \leq I_{2,j}, &&  j \in \set{J}  \\
	& 0 \leq o_{1} \leq O_{1} && \\
	& x_{1j} ,i^+_{1j}, i^-_{1j} \geq 0, && j \in \set{J}  \\
	& y_{1j} \in\{0,1\}, && j \in \set{J}. 
	\end{alignat*}
\esub
Finally, the LDDR-restricted stagewise Lagrangian dual problem is defined as
\begin{align*}
\max \  & 
\sum\limits_{t \in \set{T}} 
\Exp
\left[ 
\mathcal{L}_{t}(\alpha,\xisupt)
\right]  \\
\text{s.t.} \ & \alpha_{tk} \in \R^{J},  \qquad t \in \set{T},  k\in\set{\Kt}. 
\end{align*}


\subsection{Nonanticipative Lagrangian Dual}
Reformulate the MSLot problem as follows
\bsub
	\begin{alignat*}{2}
	\min  \ \ &   \Exp\Bigg[\sum\limits_{t \in \set{T}} \bigg[ \sum\limits_{j \in \set{J}} \Big[
	C_{tj}^{i^+}(\xisupt)i^{+,na}_{tj}(\xiT) + C_{tj}^{i^-}(\xisupt)i^{-,na}_{tj}(\xiT) + C_{tj}^{y}(\xisupt)y_{tj}^{na}(\xiT) \Big] + C_{t}^{o}(\xisupt)o_{t}^{na}(\xiT)
	\bigg]\Bigg]  \\ 
	\text{s.t.} \ \ & i^{-,na}_{tj}(\xiT)-i^{+,na}_{tj}(\xiT) + i^{+,na}_{t-1,j}(\xiT) - i^{-,na}_{t-1,j}(\xiT) + x^{na}_{t-1,j}(\xiT) = D_{tj}(\xisupt),   &&  \nonumber\\
	&&& \hspace*{-4cm} t \in \set{T}, j \in \set{J}, \allxi \\ 
	&  \sum\limits_{j \in \set{J}}(TS_jy^{na}_{tj}(\xiT) +TB_jx^{na}_{tj}(\xiT)) - o^{na}_{t}(\xiT)\leq C_t, && \hspace*{-4cm} t \in \set{T}, \allxi \\ 
	& M_{tj}\ y_{tj}^{na}(\xiT) - x_{tj}^{na}(\xiT) \geq 0, && \hspace*{-4cm} t \in \set{T}, j \in \set{J}, \allxi \\ 
	& i_{tj}^{+,na}(\xiT) \leq I_{tj}, && \hspace*{-4cm} t \in \set{T}, j \in \set{J}, \allxi \\ 
	& i_{tj}^{+,na}(\xiT) + x_{tj}^{+,na}(\xiT) \leq I_{t+1,j}, && \hspace*{-4cm} t \in \set{T}, j \in \set{J}, \allxi \\ 
	& 0 \leq o_{t}^{na}(\xiT) \leq O_{t}, &&\hspace*{-4cm}  t \in \set{T}, \allxi \\ 
	& x_{tj}^{na}(\xiT) = \Exp_{\xiprimeT}[x_{tj}^{na}(\xiprimeT)\>|\> \xisupt], && \hspace*{-4cm} t \in \set{T}, j \in \set{J}, \allxi \\ 
	& x_{tj}^{na}(\xiT), i^{+,na}_{tj}(\xiT), i^{-,na}_{tj}(\xiT) \geq 0, && \hspace*{-4cm} t \in \set{T}, j \in \set{J}, \allxi \\ 
	& y_{tj}^{na}(\xiT) \in\{0,1\}, && \hspace*{-4cm} t \in \set{T}, j \in \set{J}, \allxi  
	\end{alignat*}
\esub
where for any variable $ a^{na}(\xiT) $ the superscript $ na $ indicates the anticipative copy variable corresponding to original variable $ a(\xisupt) $.
Relaxing the nonanticipativity constraints 
using dual variables $\gamma_{t}(\xiT)$ and enforcing LDDR on these duals as
\begin{align*}
& \gamma_{tj}(\xiT) =  \sum \limits_{\kprimeindex \in \set{\Kprimet}} \Psi_{tj \kprimeindex}(\xiT) \alpha_{tj \kprimeindex},
\end{align*}
the  LDDR-restricted nonanticipative Lagrangian dual problem is obtained as follows
\begin{align*}
\max \  & 
\Exp
\left[ 
\mathcal{L}(\alpha,\xiT)
\right]  \\
\text{s.t.} \ & \alpha_{t \kprimeindex} \in \R^{J},  \qquad t \in \set{T},  \kprimeindex\in\set{\Kprimet} 
\end{align*}
where for fixed $ \hat{\alpha} $
	\begin{align*}
	\mathcal{L}(\hat{\alpha},\xiT) = \min \ \ & \sum\limits_{t \in \set{T}} \bigg[ \sum\limits_{j \in \set{J}} \Big[
	C_{tj}^{i^+}(\xisupt)i^{+,na}_{tj} + C_{tj}^{i^-}(\xisupt)i^{-,na}_{tj} + C_{tj}^{y}(\xisupt)y_{tj}^{na} \Big] + C_{t}^{o}(\xisupt)o_{t}^{na} \bigg]\nonumber\\
	+ & \sum\limits_{t \in \set{T}} \sum\limits_{j \in \set{J}} \sum \limits_{\kprimeindex \in \set{\Kprimet}} 
	\bigg( \Psi_{tj \kprimeindex}(\xiT) - \Exp_{\xiprimeT}[\Psi_{tj \kprimeindex}(\xiprimeT) | \xisupt] \bigg)
	\hat{\alpha}_{tj \kprimeindex} x_{tj}^{na} \\
	\text{s.t.} \ \   
	&  i^{-,na}_{tj} - i^{+,na}_{tj} + i^{+,na}_{t-1,j} - i^{-,na}_{t-1,j} + x^{na}_{t-1,j} = D_{tj}(\xisupt),   &&  t \in \set{T}, j \in \set{J}  \\
	&  \sum\limits_{j \in \set{J}}\big(TS_jy^{na}_{tj} +TB_jx^{na}_{tj}\big) - o^{na}_{t}\leq C_t, && t \in \set{T}   \\
	& M_{tj}\ y_{tj}^{na} - x_{tj}^{na} \geq 0, && t \in \set{T}, j \in \set{J}  \\
	& i_{tj}^{+,na} \leq I_{tj}, && t \in \set{T}, j \in \set{J}  \\
	& i_{tj}^{+,na} + x_{tj}^{+,na} \leq I_{t+1,j}, && t \in \set{T}, j \in \set{J} \\
	& 0 \leq o_{t}^{na} \leq O_{t}, && t \in \set{T}  \\
	& x_{tj}^{na}, i^{+,na}_{tj}, i^{-,na}_{tj} \geq 0, && t \in \set{T}, j \in \set{J} \\
	& y_{tj}^{na} \in\{0,1\}, && t \in \set{T}, j \in \set{J}.
\end{align*}


\subsection{Parameters of the MSLot Instances} \label{sec:appendix-parameters}
Our instances are generated using parameters that are loosely based on the work of \cite{helber2013dynamic}. Section \ref{sec:data-instances} explains the demand \revtwo{model. In this section, we specify the demand generation procedure and the rest of the parameters. Using consecutive substitutions, we can deduce that the following holds for the conditional expectation
\begin{equation*}
\Exp[D_{t+h,j} | \epsilon_{tj}, \delta_{tj}] = \mu_{t+h,j}  \left(\rho^Y\rho^{h}\left(\frac{D_{tj}-(1-\rho^Y) \, \delta_{tj}}{\rho^Y \, \mu_{tj}} - 1\right) + 1\right).
\end{equation*}

We have generated instances with $ T = 2,\dots,10 $ stages, and $ J = 3,6 $  product types. There are two sets of data instances for $ J=3 $ with  $ \rho = 0.2, \rho^Y=0.6 $ and $ \rho = 0.6, \rho^Y = 0.2 $ leading to different levels of correlation and variation among demands. The combination $ \rho = 0.6, \rho^Y = 0.2 $ is additionally used in the generation of instances with $ J = 6 $.  \rev{For solving the SW and NA dual problems, $ \lceil\frac{50}{T}\rceil |\beta|_{\#} $ and $ \lceil\frac{100}{T}\rceil |\alpha|_{\#} $ scenarios are generated respectively, and the solutions are evaluated using $ \lceil\frac{250}{T}\rceil |\alpha|_{\#} $ scenarios, where $|\cdot|_{\#}$ operator denotes the total number of variables of a given type (which depends on the choice of basis functions, discussed in Section \ref{sec:bf-selection}).}} We do not consider a production cost. The overtime cost
is $ 100 $ per unit of overtime. Holding cost is $ 15 $ per unit, while the backlog cost is $(\delta_{i^-}) (c^{i^+}_{tj}) $, where $ \delta_{i^-} = 2 $. For the last stage though, we have an end of horizon effect and $ c^{i^-}_{Tj} = 150 $. $ TB_j $ is set to $ 1 $, and 
$ TS_j = ts^{rel} \ \overline{\Exp [D_j]} \ TB_j $, where $ \overline{\Exp [D_j]} = \frac{\sum_{t\in\set{T}}\Exp[D_{tj}(\xisupt)]}{T}$ is the average expected demand of product $ j $, and $ ts^{rel} $ is  $0.25$.
The setup cost is $ c^{y}_{tj} =
\delta_y \  \overline{\Exp [D_j]} \ TBO^2 \ c^{i^+}_{ij}   $, where  $ TBO $ is the processing time between orders (set to $ 2 $), and $ \delta_y $ is $ 1.2 $.
Production capacity is $ C_t = 0.9 \frac{\sum_{j\in\set{J}}\Exp[D_{tj}(\xisupt)]}{Util}$, where $ Util $ is $0.6 $.
Inventory capacity is $ I_{tj} = \delta_{I} \  \overline{\Exp [D_j]} $, with $ \delta_{I} = 10 $.
The bound on overtime is $ O_t = \delta_{O}  C_t$, where $ \delta_{O} = 0.25 $. For big-$ M $ values in the MSLot formulation, we have
$M_j = 6 \  \overline{\Exp [D_j]} $.

\revtwo{
\subsection{NA Variable Selection}
In NA problem reformulation, another possibility to reduce the restricted dual problem size is to ignore altogether the NA constraints \eqref{eq:const-reform-nonanticipative} on some sets of variables.  For example, in order to obtain a valid formulation it is sufficient to enforce \eqref{eq:const-reform-nonanticipative} only on {\it state variables}, i.e., variables that appear in more than one decision stage, and thus one may choose to only penalize the NA constraints associated with those variables. On the other hand, in the MSLot problem, once  the recourse variables $x_t(\xisupt)$ are determined, the optimal values of the state variables $i^-(\xisupt)$ and $i^+(\xisupt)$ are immediate from the equations \eqref{eq:lot-sizing-state-equations}, suggesting it may be beneficial to penalize violation of the NA constraints on these variables. (A similar approach is used by \cite{lulli2004branch}.) A natural question for our LDDR approach is, can enforcing nonanticipativity additionally on the other variables (which is redundant before the relaxation) improve the  bound? 

In Table \ref{tab:na-var-selection}, we consider three options \rev{for instances with $J=3, \rho = 0.6, \rho^Y = 0.2$}: having either $ x $ or state variables $ i^+, i^- $ as the NA variables, or considering a combination of them. Results of Table \ref{tab:na-var-selection} show that the combination of $ x, i^+, i^- $ obtains the best bounds in a longer running time. The difference between the means of the options was tested using a $ t $-test. We find that the bounds obtained by penalizing the NA constraints on only the $ x $ variables are not significantly worse than those obtained by penalizing the NA constraints on all variables. In addition, including a penalty of the NA constraints $x$ variables leads to an improvement over just including a penalty of the NA constraints on the state variables $i^+$ and $i^-$.  Since penalizing only the NA constraints on the $x$ variables yields lower bounds that are indistinguishable from the best, in less time, we use this option in the remaining experiments.
}

\subsection{Detailed Results} \label{sec:appendix-results}
In the next pages, the actual numbers used in Section \ref{sec:results} are provided, without any scaling and normalization. Also, the confidence intervals are given for the bounds, in the form of (mean $\pm$ width).

\begin{table}[htbp]
\renewcommand{\arraystretch}{0.55}
  \centering
  \caption{Basis function selection analysis for the SW dual}
  \footnotesize{
    \begin{tabular}{crrrrrrrrr}
    \toprule
    & \multicolumn{3}{c}{\# Basis functions} & \multicolumn{3}{c}{Time (s)} & \multicolumn{3}{c}{Bound} \\
    \cmidrule(rl){2-4}\cmidrule(rl){5-7} \cmidrule(rl){8-10}
Option         & \multicolumn{1}{c}{$T=4$} & \multicolumn{1}{c}{$T=6$} & \multicolumn{1}{c}{$T=8$} & \multicolumn{1}{c}{$T=4$} & \multicolumn{1}{c}{$T=6$} & \multicolumn{1}{c}{$T=8$} & \multicolumn{1}{c}{$T=4$} & \multicolumn{1}{c}{$T=6$} & \multicolumn{1}{c}{$T=8$} \\
    \cmidrule(rl){1-1}\cmidrule(rl){2-4}\cmidrule(rl){5-7} \cmidrule(l){8-10}
    1     & 315   & 750   & 1365  & 353.8 & 1412.5 & 6273.9 & 32474.6 $\pm$ 448.6 & 41291.3 $\pm$ 1194.2 & 69849.2 $\pm$ 5151.9 \\
    2     & 180   & 300   & 420   & 91.7  & 243.7 & 632.5 & 30947.9 $\pm$ 363.9 & 39835.2 $\pm$ 496.4 & 66388.5 $\pm$ 5159.7 \\
    3     & 135   & 300   & 525   & 109   & 288.7 & 795.9 & 31053.7$\pm$ 371.9 & 40195.4 $\pm$ 671.9 & 66727.2 $\pm$ 2203.9 \\
    4     & 90    & 150   & 210   & 157.2 & 795.9 & 323.6 & 29245.5 $\pm$ 396.9 & 38958.3 $\pm$ 1086.0 & 56882.1 $\pm$ 4201.2 \\
    \bottomrule
    \end{tabular}%
}
  \label{tab:ext-sw-dual-bf}%
\end{table}%
\spacing

\begin{table}[h]
\renewcommand{\arraystretch}{0.55}
  \centering
  \caption{Basis function selection analysis for the NA dual}
  \footnotesize{
    \begin{tabular}{crrrrrrrrr}
    \toprule
     & \multicolumn{3}{c}{\# Basis functions} & \multicolumn{3}{c}{Time (s)} & \multicolumn{3}{c}{Bound} \\
     \cmidrule(rl){2-4}\cmidrule(rl){5-7} \cmidrule(rl){8-10}
Option          & \multicolumn{1}{c}{$T=4$} & \multicolumn{1}{c}{$T=6$} & \multicolumn{1}{c}{$T=8$} & \multicolumn{1}{c}{$T=4$} & \multicolumn{1}{c}{$T=6$} & \multicolumn{1}{c}{$T=8$} & \multicolumn{1}{c}{$T=4$} & \multicolumn{1}{c}{$T=6$} & \multicolumn{1}{c}{$T=8$} \\
    \cmidrule(rl){1-1}\cmidrule(rl){2-4}\cmidrule(l){5-7} \cmidrule(rl){8-10}
    1     & 756   & 1800  & 3276  & 76.9  & 403.2 & 1800.8 & 53352.5 $\pm$ 1261.5 & 71767.2 $\pm$ 2515.2 & 123173.0 $\pm$ 5101.1 \\
    2     & 432   & 720   & 1008  & 44.9  & 232.3 & 770.5 & 49940.7 $\pm$ 403.5 & 67696.9 $\pm$ 501.2 & 114913.0 $\pm$ 1834.5 \\
    3     & 324   & 720   & 1260  & 61.1  & 250   & 776.8 & 53328.5 $\pm$ 1186.6 & 71352.1 $\pm$ 1539.8 & 122755.0 $\pm$ 3085.8 \\
    4     & 216   & 360   & 504   & 34.2  & 219.5 & 735.3 & 49902.7 $\pm$ 407.0 & 67247.6 $\pm$ 442.2 & 114829.0 $\pm$ 1715.5 \\
    \bottomrule
    \end{tabular}%
}
  \label{tab:ext-na-dual-bf}%
\end{table}%
\spacing

\begin{table}[h]
\renewcommand{\arraystretch}{0.55}
  \centering
  \caption{NA variable selection analysis}
    \begin{tabular}{crrrrrr}
    \toprule
    & \multicolumn{3}{c}{Time (s)} & \multicolumn{3}{c}{Bound} \\  
    \cmidrule(rl){2-4}\cmidrule(rl){5-7}       
Option & \multicolumn{1}{c}{$T=4$} & \multicolumn{1}{c}{$T=6$} & \multicolumn{1}{c}{$T=8$} & \multicolumn{1}{c}{$T=4$} & \multicolumn{1}{c}{$T=6$} & \multicolumn{1}{c}{$T=8$} \\
    \cmidrule(rl){1-1}\cmidrule(rl){2-4}\cmidrule(rl){5-7} 
    $x$   & 61.1  & 68.1  & 175.4 & 53328.5 $\pm$ 1186.6 & 71352.1 $\pm$ 1539.8 & 122755.0 $\pm$ 3085.8 \\
    $i^+, i^-$ & 250.0   & 188.0   & 849.9 & 49456.0 $\pm$ 403.9 & 65348.7 $\pm$ 392.2 & 111605.0 $\pm$ 1778.4 \\
    $x, i^+, i^-$ & 776.8 & 653.1 & 1887.3 & 53459.1 $\pm$ 1012.2 & 71844.8 $\pm$ 939.4 & 123403.0 $\pm$ 2400.2 \\
    \bottomrule
    \end{tabular}%
  \label{tab:ext-na-select}%
\end{table}%
\spacing

\begin{table}[htbp]
\renewcommand{\arraystretch}{0.55}
  \centering
  \caption{\rev{Analysis of different values for $\lambda$ in SW UB for the \revtwo{lot-sizing} problem}}
    \begin{tabular}{crrr}
    \toprule
          & \multicolumn{3}{c}{Bound} \\
 \cmidrule(rl){2-4}$\lambda$ & \multicolumn{1}{c}{$T=4$} & \multicolumn{1}{c}{$T=6$} & \multicolumn{1}{c}{$T=8$} \\
    \cmidrule(rl){1-1} \cmidrule(rl){2-4}
    0.00  & 58292.3 $\pm$ 481.1 & 81073.9 $\pm$ 672.5 & 143503.0 $\pm$ 1912.6 \\
    0.25  & 56624.5 $\pm$ 503.2 & 80815.3 $\pm$ 689.7 & 143503.2 $\pm$ 1928.4 \\
    0.50  & 56925.1 $\pm$ 499.4 & 81138.1 $\pm$ 680.5 & 144290.8 $\pm$ 1891.7 \\
    1.00  & 118093.6 $\pm$ 1357.2 & 175553.3 $\pm$ 1921.6 & 282653.2 $\pm$ 2651.1 \\
    \bottomrule
    \end{tabular}%
  \label{tab:diff-lambda}%
\end{table}%

\newpage
\begin{table}[h]
\renewcommand{\arraystretch}{0.55}
  \centering
  \caption{Bound comparison for various instances}
    \begin{tabular}{ccrrrrrr}
    \toprule
     &  & \multicolumn{3}{c}{LB} & \multicolumn{3}{c}{UB} \\
     \cmidrule(rl){3-5} \cmidrule(rl){6-8} 
Instance          &   $ T $    & \multicolumn{1}{c}{PI} & \multicolumn{1}{c}{SW} & \multicolumn{1}{c}{NA} & \multicolumn{1}{c}{Cond Exp} & \multicolumn{1}{c}{SW} & \multicolumn{1}{c}{NA} \\
    \cmidrule(rl){1-1}\cmidrule(rl){2-2}\cmidrule(rl){3-5} \cmidrule(rl){6-8}
    \multicolumn{1}{c}{\multirow{9}[2]{*}{\begin{sideways}$\rho = 0.2, \rho^Y=0.6,J = 3$\end{sideways}}} & 2     & 46584.6 & 33327.5 & 47189.4 & 49406.4 & 49406.4 & 49001.1 \\
          & 3     & 45487.7 & 29308.6 & 46735.3 & 51076.7 & 50660.4 & 50953.9 \\
          & 4     & 48868.7 & 29962.1 & 50537.6 & 56731.2 & 54878.6 & 53614.3 \\
          & 5     & 55920.5 & 31638.4 & 56965.3 & 64759.2 & 64062.4 & 61419.5 \\
          & 6     & 64844.7 & 38036.7 & 67777.5 & 75596.2 & 75515.2 & 75384.4 \\
          & 7     & 72980.9 & 42580.4 & 75263.3 & 84230.4 & 83838.4 & $>$6h \\
          & 8     & 104296.0 & 67443.0 & 112093.0 & 129166.0 & 127026.0 & $>$6h \\
          & 9     & 102604.0 & 60101.1 & 104760.0 & 121131.0 & 117865.0 & $>$6h \\
          & 10    & 114903.0 & 70844.6 & 122374.0 & 135756.0 & 134728.0 & $>$6h \\
    \midrule
    \multicolumn{1}{c}{\multirow{9}[2]{*}{\begin{sideways}$\rho = 0.6, \rho^Y=0.2, J = 3$\end{sideways}}} & 2     & 46587.0 & 35050.1 & 46876.1 & 48484.6 & 48484.6 & 48277.4 \\
          & 3     & 45547.9 & 29282.9 & 46585.9 & 50607.9 & 50109.8 & 49404.2 \\
          & 4     & 49070.8 & 32474.3 & 53328.5 & 58292.3 & 56624.5 & 54892.8 \\
          & 5     & 55927.6 & 35362.3 & 59325.2 & 67492.5 & 66928.6 & 63767.0 \\
          & 6     & 65167.3 & 41291.3 & 71352.1 & 81073.9 & 80815.3 & 82093.8 \\
          & 7     & 73532.0 & 46230.9 & 80838.4 & 91635.6 & 90969.8 & $>$6h \\
          & 8     & 110344.0 & 69849.2 & 122755.0 & 145747.0 & 143503.0 & $>$6h \\
          & 9     & 105611.0 & 63981.7 & 115548.0 & 135670.0 & 130932.0 & $>$6h \\
          & 10    & 120331.0 & 67950.0 & 134594.0 & 155783.0 & 154407.0 & $>$6h \\
    \midrule
    \multicolumn{1}{c}{\multirow{9}[2]{*}{\begin{sideways}$\rho = 0.6, \rho^Y=0.2, J = 6$\end{sideways}}} & 2     & 97422.8 & 72505.3 & 98890.1 & 103233.0 & 103233.0 & 102963.5 \\
          & 3     & 90470.1 & 56340.1 & 92717.4 & 104545.0 & 104218.0 & 103844.1 \\
          & 4     & 94000.1 & 53199.8 & 100338 & 113226.0 & 111675.8 & 111718.4 \\
          & 5     & 107992.0 & 62849.1 & 114973.0 & 134017.0 & 132604.4 & 132614.5 \\
          & 6     & 121141.0 & 67740.8 & 131226.0 & 153875.0 & 153151.2 & $>$6h \\
          & 7     & 118421.0 & 65444.7 & 130542.0 & 175972.0 & 174715.2 & $>$6h \\
          & 8     & 139644.0 & 106919.0 & 157804.0 & 240137.0 & 238465.8 & $>$6h \\
          & 9     & 127301.0 & 107217.7 & 142417.0 & 235657.0 & 233470.2 & $>$6h \\
          & 10    & 135019.0 & 118333.6 & 156096.0 & 273716.0 & 270377.2 & $>$6h \\
    \bottomrule
    \end{tabular}%
  \label{tab:ext-bound-comp}%
\end{table}%
\spacing

\newpage
\begin{table}[h]
\renewcommand{\arraystretch}{0.55}
  \centering
  \caption{\rev{Solution time of bounding techniques in seconds for the MSLot problem. For UB techniques, the average time of one sample path is given. For NA upper bounds with $T>6$, the reported value is the average time over an evaluation sample of size 10.}}
\begin{tabular}{rcrrrr}
\toprule
      &       & \multicolumn{2}{c}{SW Time (s)} & \multicolumn{2}{c}{NA  Time (s)} \\
\cmidrule(rl){3-4}\cmidrule(rl){5-6}\multicolumn{1}{c}{Instance} & $T$   & \multicolumn{1}{c}{LB} & \multicolumn{1}{c}{UB (avg./sample path)} & \multicolumn{1}{c}{LB} & \multicolumn{1}{c}{UB (avg./sample path)} \\
\cmidrule(rl){1-1}\cmidrule(rl){2-2}\cmidrule(rl){3-4}\cmidrule(rl){5-6}
       \multicolumn{1}{c}{\multirow{9}[2]{*}{\begin{sideways}$\rho = 0.2, \rho^Y=0.6,J = 3$\end{sideways}}} & 2     & 3.4   & $<$0.1s & 3.9   & $<$0.1s \\
      & 3     & 18.3  & $<$0.1s & 24.4  & 0.2 \\
      & 4     & 161.4 & $<$0.1s & 59.7  & 0.8 \\
      & 5     & 338.3 & $<$0.1s & 115.9 & 1.7 \\
      & 6     & 750.7 & 0.1   & 339.0 & 3.7 \\
      & 7     & 960.0 & 0.2   & 437.4 & 7.7 \\
      & 8     & 2184.1 & 0.3   & 605.3 & 9.5 \\
      & 9     & 2460.0 & 0.4   & 714.0 & 12.4 \\
      & 10    & 2711.4 & 0.5   & 1975.4 & 17.6 \\
\midrule
      \multicolumn{1}{c}{\multirow{9}[2]{*}{\begin{sideways}$\rho = 0.6, \rho^Y=0.2, J = 3$\end{sideways}}} & 2     & 3.5   & $<$0.1s & 3.9   & $<$0.1s \\
      & 3     & 17.6  & $<$0.1s & 25.2  & 0.2 \\
      & 4     & 67.6  & $<$0.1s & 58.2  & 0.7 \\
      & 5     & 115.1 & $<$0.1s & 72.5  & 1.8 \\
      & 6     & 675.2 & 0.1   & 374.6 & 4.0 \\
      & 7     & 883.8 & 0.1   & 413.5 & 7.4 \\
      & 8     & 1499.0 & 0.3   & 645.0 & 9.6 \\
      & 9     & 1725.6 & 0.5   & 696.5 & 12.2 \\
      & 10    & 1919.7 & 0.4   & 1741.7 & 16.5 \\
\midrule
      \multicolumn{1}{c}{\multirow{9}[2]{*}{\begin{sideways}$\rho = 0.6, \rho^Y=0.2, J = 6$\end{sideways}}} & 2     & 20.2  & $<$0.1s & 19.0  & 0.1 \\
      & 3     & 288.6 & $<$0.1s & 119.8 & 0.6 \\
      & 4     & 512.7 & $<$0.1s & 679.9 & 2.0 \\
      & 5     & 4747.9 & $<$0.1s & 1334.4 & 5.0 \\
      & 6     & 6154.1 & 0.1   & 2290.8 & 8.2 \\
      & 7     & 10963.4 & 0.2   & 3749.9 & 14.2 \\
      & 8     & 12859.1 & 0.3   & 7805.3 & 25.6 \\
      & 9     & 13946.1 & 0.5   & 8912.0 & 29.9 \\
      & 10    & 15524.3 & 0.5   & 11228.1 & 42.0 \\
\bottomrule
\end{tabular}%
  \label{tab:time_lotsizing}%
\end{table}%


\section{Solution Times for MSRWA Instances}
\begin{table}[h]
\renewcommand{\arraystretch}{0.55}
  \centering
  \caption{\revtwo{Solution time of bounding techniques in seconds for the MSRWA problem.}}
    \begin{tabular}{lcrcrc}
    \toprule
          &       & \multicolumn{2}{c}{Basis Function Option 1 Time (s)} & \multicolumn{2}{c}{Basis Function Option 2 Time (s)} \\
\cmidrule(lr){3-4} \cmidrule(lr){5-6}    $(|\waveSet|,\bar{\traffic},\gamma_0,\delta_1)$ & $T$   & \multicolumn{1}{c}{NA UB} & \multicolumn{1}{c}{SW LB (avg./sample path)} & \multicolumn{1}{c}{NA UB } & \multicolumn{1}{c}{SW LB (avg./sample path)} \\
    \cmidrule(lr){1-1} \cmidrule(lr){2-2}
    \cmidrule(lr){3-4} \cmidrule(lr){5-6}
    \multirow{3}[-2]{*}{$(1,10,3,0.2)$} & 3     & 122.2 & 0.2   & 76.0  & 0.2 \\
          & 4     & 417.0 & 0.4   & 200.2 & 0.4 \\
          & 5     & 1852.0 & 0.8   & 1193.2 & 0.9 \\
    \midrule
    \multirow{3}[-2]{*}{$(1,10,4.5,0.05)$} & 3     & 113.5 & 0.2   & 93.1  & 0.2 \\
          & 4     & 420.4 & 0.4   & 277.6 & 0.5 \\
          & 5     & 1325.7 & 0.8   & 951.0 & 0.9 \\
    \midrule
    \multirow{3}[-2]{*}{$(3,30,9,0.2)$} & 3     & 772.2 & 0.2   & 302.3 & 0.3 \\
          & 4     & 1533.4 & 0.5   & 787.0 & 0.5 \\
          & 5     & 5240.1 & 0.9   & 3339.9 & 0.9 \\
    \midrule
    \multirow{3}[-2]{*}{$(3,30,13.5,0.05)$} & 3     & 945.4 & 0.2   & 494.6 & 0.2 \\
          & 4     & 2800.8 & 0.5   & 1045.2 & 0.5 \\
          & 5     & 7117.1 & 0.9   & 3061.1 & 0.9 \\
    \bottomrule
    \end{tabular}%
  \label{tab:time_rwa}%
\end{table}%

\end{appendices}

\end{document}